\documentclass[reqno,10pt, centertags,draft]{amsart}
\usepackage{amsmath,amsthm,amscd,amssymb,latexsym,upref}

\makeatletter
\def\theequation{\@arabic\c@equation}

\newcommand{\bbN}{{\mathbb{N}}}
\newcommand{\bbR}{{\mathbb{R}}}

\newcommand{\bbZ}{{\mathbb{Z}}}
\newcommand{\bbC}{{\mathbb{C}}}

\newcommand{\cA}{{\mathcal A}}

\newcommand{\cH}{{\mathcal H}}
\newcommand{\cI}{{\mathcal I}}

\newcommand{\cM}{{\mathcal M}}

\newcommand{\cP}{{\mathcal P}}
\newcommand{\cR}{{\mathcal R}}
\newcommand{\cS}{{\mathcal S}}
\newcommand{\cT}{{\mathcal T}}

\newcommand{\no}{\nonumber}
\newcommand{\lb}{\label}
\newcommand{\f}{\frac}

\newcommand{\spec}{\text{\rm{spec}}}

\newcommand{\ran}{\text{\rm{ran}}}

\newcommand{\dom}{\text{\rm{dom}}}

\newcommand{\supp}{\text{\rm{supp}}}

\newcommand{\bi}{\bibitem}

\newcommand{\Ln}{\text{\rm{Ln}}}
\newcommand{\tr}{\text{\rm{tr}}}
\newcommand{\Aut}{\text{\rm{Aut}}}
\newcommand{\Int}{\text{\rm{int}}}

\renewcommand{\Re}{\text{\rm Re}}
\renewcommand{\Im}{\text{\rm Im}}

\numberwithin{equation}{section}

\newtheorem{theorem}{Theorem}[section]
\newtheorem{lemma}[theorem]{Lemma}
\newtheorem{corollary}[theorem]{Corollary}
\newtheorem{hypothesis}[theorem]{Hypothesis}
\theoremstyle{definition}
\newtheorem{definition}[theorem]{Definition}
\theoremstyle{remark}
\newtheorem{remark}[theorem]{Remark}
\newtheorem{example}[theorem]{Example}

\begin{document}

\title[$SL_2(\bbR)$, Herglotz functions, and 
spectral averaging]{$\mathbf {SL_2(\bbR)}$, Exponential Herglotz
Representations, \\ and Spectral Averaging}
\author[F.\ Gesztesy and K.\ A.\ Makarov]{Fritz Gesztesy and Konstantin
A.\ Makarov}
\address{Department of Mathematics,
University of Missouri, Columbia, MO 65211, USA}
\email{fritz@math.missouri.edu} 
\urladdr{http://www.math.missouri.edu/people/fgesztesy.html}
\address{Department of Mathematics, University of
Missouri, Columbia, MO 65211, USA}
\email{makarov@math.missouri.edu} 
\urladdr{http://www.math.missouri.edu/people/kmakarov.html}
\date{April 19, 2002}
\thanks{{\it St.\ Petersburg Math.\ J.} (to appear)}

\begin{abstract}
We revisit the concept of spectral averaging and point out its origin in
connection with one-parameter subgroups of $SL_2(\bbR)$ and
the corresponding M\"obius transformations. In particular, we identify
exponential Herglotz representations as the basic ingredient for the 
absolute continuity of average spectral measures with respect to
Lebesgue measure and the associated spectral shift function as the
corresponding density for the averaged measure. As a by-product of our
investigations we unify the treatment of rank-one perturbations of
self-adjoint operators and that of self-adjoint extensions of symmetric
operators with deficiency indices $(1,1)$. Moreover, we derive separate 
averaging results for absolutely continuous, singularly continuous, and
pure point measures and conclude with an averaging result of the
$\kappa$-continuous part (with respect to the $\kappa$-dimensional
Hausdorff measure) of singularly continuous measures.
\end{abstract}

\maketitle

\section{Introduction} \lb{s1}

Spectral averaging is usually concerned with integrating the spectral
measure of a  one-parameter family of self-adjoint operators with
respect
to a parameter, typically a coupling constant or a boundary condition
parameter. One then proceeds to proving the absolute continuity of
the  integrated (averaged) spectral measure with respect to Lebesgue
measure. Actually, one is usually more ambitious and tries to establish 
the universality of spectral averaging, provided that averaging is
carried
out over the whole parameter space. That is, one intends to prove that
the
averaged measure does not depend upon the concrete choice of the
one-parameter family of operators and that it is mutually equivalent to
Lebesgue measure.

In this paper we revisit this circle of ideas and present a discussion
of the following topics: 

\smallskip

\noindent $\bullet$ The intimate connection between spectral 
averaging, $SL_2(\bbR)$, and M\"obius transformations. \\
$\bullet$ The exponential Herglotz representation theorem is shown
to be the underlying reason for absolute continuity of averaged
spectral measures with respect to Lebesgue measure. In particular,
this identifies the spectral shift function as the density in the
absolutely continuous averaged spectral measure. \\
$\bullet$ Various existing results on (universality of)
spectral averaging are extended. In particular, we don't assume the
existence of a spectral gap (or boundedness from below) in the
associated self-adjoint operators. \\
$\bullet$ Conditions for (non)universality of spectral
averaging to hold are identified. \\
$\bullet$ A unified treatment of self-adjoint rank-one
perturbations of a self-adjoint operator and self-adjoint 
extensions of a densely defined closed symmetric operator with
deficiency indices $(1,1)$ is presented. \\
$\bullet$ Separate averaging for point spectra, absolutely continuous, 
and singularly continuous spectra are discussed. \\
$\bullet$ A partial result for averaging the $\kappa$-continuous
part (with respect to the $\kappa$-dimensional Hausdorff measure)
of singularly continuous measures is derived. 

\medskip

We next illustrate these ideas in two canonical cases: the case of 
rank-one perturbation theory and that of the theory of self-adjoint
extensions of symmetric operators  with deficiency indices $(1,1)$.

Let $A$ be a self-adjoint operator in a separable complex Hilbert space 
$\cH$ and $P$ an orthogonal rank-one projection in $\cH$. We introduce
two Herglotz functions, $M$ and $N$, associated with the pair $(A,P)$
\begin{equation}\label{MM}
M(z)=\tr(P(A-z)^{-1}P), \quad z\in\bbC_+
\end{equation}
and 
\begin{equation}\label{NN}
N(z)=\tr(P(zA+I)(A-z)^{-1}P), \quad z\in\bbC_+,
\end{equation}
with $\bbC_+$ the open upper complex half-plane. One then has the
Herglotz
representations, 
\begin{equation}
M(z)=\int_\bbR\frac{d\mu(\lambda)}{\lambda-z}, \quad z\in\bbC_+, \lb{mm}
\end{equation}
with $\mu$ a probability measure on $\bbR$, $\mu(\bbR)=1$, and 
\begin{equation}
N(z)=B+\int_\bbR d\nu(\lambda) \bigg (\frac{1}{\lambda-z}
-\frac{\lambda}{1+\lambda^2}\bigg ), \quad z\in\bbC_+, \lb{Nn}
\end{equation}
with $B\in\bbR$ and $\nu$ a Borel measure satisfying 
\begin{equation}
\int_\bbR \frac{d\nu(\lambda)}{1+\lambda^2}<\infty. \lb{nu}
\end{equation}
Actually, a short computation reveals that 
\begin{equation}
N(z)= z +(1+z^2)M(z) 
\end{equation}
and 
\begin{equation}
d\nu(\lambda)=(1+\lambda^2)d\mu(\lambda), \quad B=\Re(N(i))=0.
\end{equation}
Thus, \eqref{Nn} simplifies to
\begin{equation}
N(z)=\int_\bbR (1+\lambda^2)d\mu(\lambda) \bigg (\frac{1}{\lambda-z}
-\frac{\lambda}{1+\lambda^2}\bigg ), \quad z\in\bbC_+. \lb{nn}
\end{equation}
If $A$ is unbounded and $\ran(P)\cap\dom(A)=\{0\}$, then the measure
$\nu$
is infinite.

\begin{lemma} \lb{l1.1}
Consider the one-parameter family of self-adjoint operators, 
\begin{equation}
A_t=A+t P, \quad t\in\bbR, 
\end{equation}
with resolvents 
\begin{equation}
(A_t-z)^{-1}=(A-z)^{-1}-\frac{1}{ M(z)+(1/t)}(A-z)^{-1}P(A-
z)^{-1}, \quad (t,z)\in\bbR\times\bbC_+,
\end{equation}
and $M$ given by \eqref{MM}. Moreover, introduce
\begin{equation}\label{pert}
M_t(z)=\frac{M(z)}{t M(z)+1}, \quad (t,z)\in\bbR\times\bbC_+.
\end{equation}
Then $M_t$ is the corresponding $M$-function associated with the pair
$(A_t, P)$ $($cf.\ \eqref{MM}$)$. Denote by $\mu_t$ the measure in
\eqref{mm} associated with $M_t$ and by $\Delta$ a bounded Borel set on
$\bbR$. Then averaging $\mu_t$ yields an absolutely continuous measure
with respect to Lebesgue measure,  
\begin{equation}\label{BS}
\int_{t_1}^{t_2}d t \, \mu_t (\Delta)=\int_\Delta
d\lambda \, [\xi(\lambda;A_{t_2}, A)-\xi(\lambda;A_{t_1}, A)], 
\end{equation}
where $\xi(\cdot,B,A)$ is the spectral shift function associated
with the pair $(B,A)$ of self-adjoint operators. Moreover, spectral
averaging is universal in the sense that 
\begin{equation}
\int_{-\infty}^\infty dt \, \mu_t (\Delta)=|\Delta|, \lb{up}
\end{equation}
with $|\cdot|$ denoting Lebesgue measure on $\bbR$.
\end{lemma}

\begin{remark}
The proof of the lemma is well-known and can be found in \cite{GS95}
and
\cite{Si95}. In fact, \eqref{BS} is a particular case of the
Birman--Solomyak spectral averaging  formula \cite{BS75} proven in the
mid-seventies.
\end{remark}

\begin{lemma} \lb{l1.3}
Assume that $A$ is an unbounded self-adjoint operator and
\begin{equation}\label{RD}
\ran(P)\cap \dom(A)=0.
\end{equation}
Then the one-parameter family of operator-valued functions 
\begin{align}\label{Krein}
R_t(z)&=(A-z)^{-1}
-\frac{1}{N(z)+(1/t)}
(A-i)(A-z)^{-1}P(A+i)(A-z)^{-1}, \\ 
& \hspace*{7.6cm} (t,z)\in\bbR\times\bbC_+, \no 
\end{align}
with $N$ given by \eqref{NN}, are the resolvents of a one-parameter
family of self-adjoint operators $\{A_t\}_{t\in\bbR}$. Introducing,
\begin{equation}\label{ext}
N_t(z)=\frac{N(z)-t  }{t N(z)+1}, \quad (t,z)\in\bbR\times\bbC_+,
\end{equation}
then $N_t$ is the $N$-function $($in the sense of \eqref{Krein}$)$ of
the pair $(A_t, P)$ $($cf.\ \eqref{nn}$)$. The family
$\{A_t\}_{t\in\bbR}$ is a one-parameter family of self-adjoint
extensions of a closed symmetric densely defined operator $\dot A$
with deficiency indices $(1,1)$, 
\begin{equation}
\dot A=A\big|_{\dom(\dot A)}, \quad 
\dom(\dot A)=\bigcap_{t\in \bbR}\dom(A_t). 
\end{equation}
In particular, $\lim_{t\to 0}A_{t}=A$ in the strong resolvent sense.
Denote
by $\nu_t$ the measure in \eqref{nn} associated with $N_t$ and by
$\Delta$
a bounded Borel set on $\bbR$. Then averaging $\nu_t$ yields an
absolutely continuous measure with respect to Lebesgue measure,
\begin{equation}\label{Ja}
\frac{1}{\pi}\int_{t_1}^{t_2} \frac{dt}{1+t^2} \,
\nu_t(\Delta)=\int_\Delta
d\lambda \, [\xi(\lambda;A_{t_2}, A)-\xi(\lambda;A_{t_1}, A)].
\end{equation}
Moreover, spectral averaging is universal in the sense that 
\begin{equation}
\frac{1}{\pi}\int_{-\infty}^{\infty} \frac{dt}{1+t^2} \, \nu_t
(\Delta)=|\Delta|. \lb{ue}
\end{equation}
\end{lemma}

\begin{remark}
The resolvent formula \eqref{Krein} is due to Krein \cite{Kr44} and
Naimark \cite{Na43}. The proof of the  transformation law \eqref{ext}
can be found, for instance, in \cite{Do65}. The spectral averaging
formula \eqref{Ja} in the case of boundary condition dependence for a
semibounded Schr\"odinger operator is due to Javrjan \cite{Ja71}.
Javrjan's method can easily be adapted to the case of arbitrary
self-adjoint operators $A$ having  a spectral gap. The treatment of
the general case of $A$ with $\spec(A)=\bbR$ needs some additional
information on the spectral shift function theory in the case of 
relatively trace class perturbation. In this case the  spectral shift
function should be viewed  as a   path-dependent homotopy invariant
characteristics of the perturbation $($see,
\cite[Ch.\ 8, Sect.\ 8]{Ya92}$)$ and the proof of \eqref{Ja}
requires minor additional efforts.
\end{remark}

In the case of perturbation theory   the transformation \eqref{pert} can
be represented in the form
\begin{equation}
M_t(z)=g_t(M(z)), 
\end{equation}
where $\{g_t\}_{t\in \bbR}$ is a one-parameter group of automorphisms of
the open upper half-plane $\bbC_+$
\begin{equation}\label{struc}
g_t\circ g_s=g_{t+s}, \quad s,t\in\bbR,
\end{equation}
where
\begin{equation}
g_t(z)=\frac{z}{t z+1}, \quad (t,z)\in \bbR\times\bbC_+. 
\end{equation}
In the case of self-adjoint extension theory the transformation
\eqref{ext}
can be written as
\begin{equation}
N_t(z)=f_t(N(z)), 
\end{equation}
where $\{f_t\}_{t\in \bbR}$ is a one-parameter family of
automorphisms of $\bbC_+$
\begin{equation}
f_t(z)=\frac{z-t}{t z+1}, \quad (t,z)\in \bbR\times\bbC_+. 
\end{equation} 
The family of transformations $\{f_t\}_{t\in \bbR}$ is not a
one-parameter
subgroup of $SL_2(\bbR)$. However, by a change of parametrization $t
\mapsto \tan (t)$, the group law \eqref{struc} can be restored with 
\begin{equation}
g_t(z)=f_{\text{arctan}(t)}(z),  \quad (t,z)\in \bbR\times\bbC_+. 
\lb{1.23} 
\end{equation}
In either case, the one-parameter family $\{g_t\}_{t\in\bbR}$ of
automorphisms of $\bbC_+$ gives rise to a dynamical system on a
certain ``phase space'' of measures as discussed in Section \ref{s3}.  

We continue with an intuitive explanation of how exponential Herglotz
representations, and hence spectral shift functions, naturally enter 
the averaging process \eqref{BS}, \eqref{Ja}. In both cases, Lemma
\ref{l1.1} and \ref{l1.3}, $M_t$, respectively, $N_t$ (the latter after
reparametrizing $t\mapsto \tan(t)$) are of the type,
\begin{equation}
M_t(z)=\f{a_t M_0(z) +b_t}{c_t M_0(z) + d_t}
=\f{d}{dt}\Ln(c_t M_0(z)+d_t), \quad (t,z)\in \bbR\times\bbC_+. 
\end{equation}
Here $M_0$ represents $M$ and $N$ in Lemmas \ref{l1.1} and \ref{l1.3},
respectively, $\Ln(\cdot)$ denotes the logarithm on the standard
infinitely
sheeted Riemann surface branched at zero and infinity (and some care
taking appropriate sheets must be exercised), and the coefficients
$a_t$, $b_t$, $c_t$, $d_t$ are all real-valued satisfying
\begin{equation}
\begin{pmatrix} a_t & b_t \\ c_t & d_t \end{pmatrix}\bigg|_{t=0}
=\begin{pmatrix} 1 & 0 \\ 0 & 1 \end{pmatrix} \text{ and } 
a_td_t-b_tc_t=1, \; t\in\bbR.
\end{equation}
Since $M_0$ is a Herglotz function, so is $M_t$ for each $t\in\bbR$.
Similarly, $c_t M_0 +d_t$ is a Herglotz or anti-Herglotz function and
thus
$M_t$ and $c_t M_0 +d_t$ admit Herglotz and exponential Herglotz
representations of the type,
\begin{align}
M_t(z)&=B_t+\int_\bbR d\, \omega_t(\lambda)\bigg
(\frac{1}{\lambda-z}-\frac{\lambda}{1+\lambda^2}\bigg ), \\
\Ln(c_tM_0(z)+d_t))&=C_t+\int_\bbR d\lambda \, \xi_t(\lambda)\bigg
(\frac{1}{\lambda-z}-\frac{\lambda}{1+\lambda^2}\bigg ), 
\end{align}
where $B_t, C_t\in\bbR$,
\begin{align}
\omega_t((\lambda_1,\lambda_2])&=\lim_{\delta\downarrow 0}
\lim_{\varepsilon\downarrow 0} \f{1}{\pi}
\int_{\lambda_1+\delta}^{\lambda_2+\delta}
d\lambda\,\Im(M_t(\lambda+i\varepsilon)), \lb{omega} \\
\xi_t(\lambda)&=\frac{1}{\pi}\lim_{\varepsilon \downarrow 0}
\Im\big (\Ln(c_tM_0(\lambda+i \varepsilon)+d_t)\big) 
\text{ for a.e.\ $\lambda\in\bbR$,} \lb{xI} 
\end{align}
and 
\begin{equation}
\int_{\bbR} \f{d\omega_t (\lambda)}{1+\lambda^2}<\infty, \quad
\xi_t(\cdot)\in L^\infty(\bbR), \quad t\in\bbR.
\end{equation}
Thus, one formally obtains for any bounded Borel set
$\Delta\subset\bbR$, 
\begin{align}
\int_{t_1}^{t_2} dt\,\omega_t(\Delta)&=\f{1}{\pi}\int_{\Delta}d\lambda 
\int_{t_1}^{t_2} dt\, \lim_{\varepsilon\downarrow 0} \f{d}{dt} 
\Im\big(\Ln(c_t M_0(\lambda+i\varepsilon)+d_t)\big) \no \\
&= \int_{\Delta} d\lambda \int_{t_1}^{t_2} dt\, \f{d}{dt}
\xi_t(\lambda) \no \\
&=\int_{\Delta} d\lambda \,
[\xi_{t_2}(\lambda)-\xi_{t_1}(\lambda)], \lb{ave}
\end{align}
freely interchanging integrals, limits, and differentiation. Once
rigorously established, \eqref{ave} proves that averaging $\omega_t$
over
the interval $[t_1,t_2]$ yields a measure absolutely continuous with
respect to Lebesgue measure on $\bbR$ and density related to the
spectral
shift function $\xi=\xi_{t_2}-\xi_{t_1}$. Moreover, in the case of
perturbations discussed in Lemma \ref{l1.1}, one can show that   
\begin{equation}
\xi_{t_2}(\lambda)-\xi_{t_1}(\lambda)\to 1 \, \text{ as
$t_1\downarrow -\infty$ and $t_2\uparrow\infty$}
\end{equation}
and hence the universal behavior \eqref{up}  
\begin{equation}
\int_{-\infty}^\infty dt\, \omega_t(\Delta)=|\Delta|
\end{equation}
emerges. The case of self-adjoint extensions discussed in Lemma
\ref{l1.3} requires some additional periodicity considerations with
respect to $t$ but in the end also yields the universality in
\eqref{ue}. However, a third case of one-parameter subgroups of
$SL_2(\bbR)$ considered in the following sections shows that
universality cannot be taken for granted and may in fact fail. The
material in Sections \ref{s2} and \ref{s3} will justify the formal
procedures in \eqref{ave}. 

Before describing the contents of each section we briefly review 
the historical development of this subject, which appears to be
less well-known. To the best of our knowledge, the credit for the
first paper on spectral averaging belongs to Javrjan \cite{Ja66} (see
also the subsequent \cite{Ja71}), who
considered half-line Schr\"odinger operators on $(0,\infty)$ and
averaged over the boundary condition parameter at $x=0$ as early as
1966. The next step is  due to Birman and Solomyak \cite{BS75} in
1975. They considered trace class perturbations of self-adjoint
operators and averaged over the coupling constant parameter (using 
the differentiation formula for operator-valued functions
by Dalecki{\u\i} and S.~Kre\u\i n \cite{DK51}). Aleksandrov
\cite{Al87} appears to be the first to consider spectral averaging of
a measure and separately averaging of its singular part in connection
with the boundary behavior of inner functions in the unit disk in
1987. More recent treatments of spectral averaging can be found in
Birman and Pushnitski \cite{BP98}, Gesztesy and Makarov \cite{GM99},
Gesztesy Makarov, and Naboko \cite{GMN99} (the latter references
discuss an operator-valued version of the Birman--Solomyak averaging
formula), Gesztesy, Makarov, and Motovilov \cite{GMM00}, and Simon
\cite{Si95}, 
\cite{Si98}.

The concept of spectral averaging became an important tool in
investigations of disordered systems, in particular, in connection with
random Schr\"odinger and Jacobi operators in the early eighties. In
1983, Carmona \cite{Ca83} (see also \cite{Ca84}), apparently unaware of
previous results by Javrjan and Birman and Solomyak, used spectral
averaging over boundary condition parameters to prove the existence of
an absolutely continuous (a.c.) component in random and deterministic
Schr\"odinger operators (he also proved that the rest of the spectrum
consists of eigenvalues dense in certain intervals with exponentially
localized eigenfunctions in some random cases). Kotani also used
this approach to link the existence of pure point spectrum and
exponentially decaying eigenfunctions with the positivity of the
Lyapunov exponent in 1984 \cite{Ko86} (published in 1986). Kotani's
work inspired new proofs of exponential localization by Delyon, L\'evy,
and Souillard \cite{DLS85}, \cite{DLS85a}, Simon and Wolff
\cite{SW86}, Simon \cite{Si85}, Delyon, Simon and Souillard
\cite{DSS87}, and Kotani and Simon \cite{KS87} for one- and
quasi-one-dimensional as well as multi-dimensional Anderson models
(the latter for large disorder or sufficiently high energy) and
one-dimensional random Schr\"odinger operators. In all these
references spectral averaging over coupling constants plays a crucial
role. This is especially transparent in the paper by Simon and Wolff
\cite{SW86}, which uses results by Aronszajn \cite{Ar57} and Donoghue
\cite{Do65} as their point of departure to study the variation of
singular spectra under rank-one perturbations of self-adjoint
operators. This is also discussed in Simon's review \cite{Si95}. (For
textbook presentations of spectral averaging in this context we refer
to \cite[Sect.\ VIII.2]{CL90}, \cite[Sect.\ 13]{PF92}.) Subsequently,
Gordon \cite{Go94}, \cite{Go97} used spectral averaging in his studies
of eigenvalues embedded in the essential spectrum. Spectral averaging
has also been used to prove exponential localization for the
one-dimensional Poisson model by Stolz
\cite{St95}. A more general approach, involving two-parameter spectral
averaging, has recently been employed to prove exponential
localization in the Poisson and random displacement models in one
dimension by Buschmann and Stolz \cite{BS01}. The latter approach was
again used by Sims and Stolz \cite{SS00} in their discussion of
exponential localization of the one-dimensional random displacement
model and in a one-dimensional model of wave propagation in a random
medium. Combes and Hislop \cite{CH94} use averaging of spectral
families to prove a Wegner-type estimate for a family of Anderson
and Poisson-like multi-dimensional random Hamiltonians. Moreover,
Combes, Hislop, and Mourre \cite{CHM96} in their discussion of
perturbations of singular spectra and exponential localization for
certain multi-dimensional random Schr\"odinger operators, and Combes,
Hislop, Klopp, and Nakamura \cite{CHKN01} in their study of the Wegner
estimate and the integrated density of states, discuss spectral
averaging in the spirit of Birman and Solomyak.

In Section \ref{s2} we collect basic facts on $SL_2(\bbR)$, 
M\"obius transformations, and the infinitely sheeted Riemann surface of
the logarithm, as needed in the subsequent sections. Section \ref{s3},
the principal section of this paper, then develops spectral averaging
for spectral measures as well as for the associated absolutely
continuous, singularly contionuous, and pure point parts (with respect
to Lebesgue measure). Finally, Section \ref{s4} obtains a partial
result concerning spectral averaging of the $\kappa$-continuous
part (with respect to the $\kappa$-dimensional Hausdorff measure)
of the singularly continuous part of measures.

\section{Preliminaries on $SL_2(\bbR)$ and on M\"obius transformations}
\lb{s2}

$SL_2(\bbR)$ denotes the group of $2\times 2$ real matrices with
determinant equal to $1$. By definition, its Lie
algebra, $\mathit{sl}_2(\bbR)$,  consists of those matrices $X$ such that
$e^{tX}\in SL_2(\bbR)$ for all $t\in \bbR$ (cf., e.g., \cite[Ch.\
VI]{La85}). Therefore, $\mathit{sl}_2(\bbR)$ consists of all $2\times 2$
real matrices $X$ with zero trace, $\tr(X)=0$. The following three
matrices then form a basis  for $\mathit{sl}_2(\bbR)$
\begin{equation}
X_1=\begin{pmatrix} 0 & 1 \\ 0 & 0 \end{pmatrix}, \quad 
X_2=\begin{pmatrix} 1 & 0 \\ 0 & -1 \end{pmatrix}, \quad 
X_3=\begin{pmatrix} 0 & 0 \\ 1 & 0 \end{pmatrix},
\end{equation}
and one verifies the following commutation relations
\begin{equation}
[X_2, X_1]=2X_1, \quad [X_1, X_3]=X_2, \quad [X_3, X_2]=2X_3.
\end{equation}
If $X\in\mathit{sl}_2(\bbR)$ then the map
$t\mapsto e^{tX}$, $t\in \bbR$ is a one-parameter subgroup of
$SL_2(\bbR)$ and all one-parameter subgroups can be obtained in that
way.

For future reference we recall the notion of automorphisms of
the open complex upper half-plane $\bbC_+$, denoted by $\Aut(\bbC_+)$:
\begin{equation}
\Aut(\bbC_+)=\{g\colon \bbC_+\to\bbC_+\,|\, \text{$g$ is biholomorphic 
(i.e., a conformal self-map of $\bbC_+$)}\}. \lb{2.2}
\end{equation}
$\Aut(\bbC_+)$ becomes a group with respect to compositions of maps.  
For simplicity, this group is denoted by the same symbol.

To fix the notational setup we now introduce the following hypothesis.

\begin{hypothesis}\label{auto}
Given $\alpha,\beta,\gamma \in\bbR$, represent an element
$X=X(\alpha,\beta,\gamma)\in\mathit{sl}_2(\bbR)$ as 
\begin{equation}
X= \alpha X_1+\beta X_2 +\gamma X_3=
\begin{pmatrix} \beta & \alpha \\ \gamma & -\beta \end{pmatrix} \lb{2.3}
\end{equation}
and denote by 
\begin{align}
\begin{split}
g_t(z)&= \frac{a_t z+b_t}{c_t z +d_t}, \quad (t,z)\in\bbR\times\bbC_+,
\\
g_0(z)&=z, \quad z\in\bbC_+ \lb{g}
\end{split}
\end{align}
the corresponding one-parameter group of automorphisms of the open 
upper-half plane $\bbC_+$ such that 
 \begin{equation} \label{entry}
\begin{pmatrix} a_t & b_t  \\ c_t & d_t \end{pmatrix}= e^{tX}\in
SL_2(\bbR), \quad t\in \bbR.
\end{equation}
\end{hypothesis}

We briefly recall a few facts in connection with M\"obius (i.e., linear
fractional) transformations \eqref{g}. Let $M$ be a M\"obius
transformation of the type 
\begin{equation}
M(z)=\f{az+b}{cz+d}, \quad z\in\bbC\cup\{\infty\}, \; a,b,c,d\in\bbC, \,
ad-bc\neq 0. \lb{Mob}
\end{equation}
Then, \\
(i) $M$ maps $\bbR\cup\{\infty\}$ onto itself if and only if $M$ admits
a
representation where $a,b,c,d\in\bbR$ and $|ad-bc|=1$. \\
(ii) $M$ maps $\bbC_+$ onto itself if and only if $M$ admits a
representation where $a,b,c,d\in\bbR$ and $ad-bc=1$. \\
(iii) $\Aut(\bbC_+)$ is isomorphic to $SL_2(\bbR)/\{I_2,-I_2\}$
($I_2$ the identity matrix in $\bbR^2$). \\
(iv) Assuming $\det(M)=ad-bc=1$ in \eqref{Mob}, one uses $\tr(M)=(a+d)$
to
classify $M$ as 

{\it elliptic}, if $(a+d)\in\bbR$ and $|a+d|<2$ 

{\it parabolic}, if $(a+d)=\pm 2$ 

{\it hyperbolic}, if $(a+d)\in\bbR$ and $|a+d|>2$ 

{\it loxodromic}, if $(a+d)\in\bbC\backslash\bbR$.
\medskip
  
\noindent On the other hand, assuming $\left(\begin{smallmatrix} a_t &
b_t
\\ c_t & d_t \end{smallmatrix}\right)=e^{tX}$, $t\in\bbR$, with
$\tr(X)=0$
and
$X=\left(\begin{smallmatrix} \beta & \alpha \\ \gamma & -\beta
\end{smallmatrix}\right)$, one can use
$\det(X)=-\alpha\gamma-\beta^2$ to classify the one-parameter subgroups
of M\"obius transformations in \eqref{g} and distinguish three cases: 

{\it Case} I: $\,\,\,\,\,\det (X)>0$ (cyclic subgroup)

{\it Case} II: $\,\,\,\det (X)=0$ 

{\it Case} III: $\,\det (X)<0$ (hyperbolic subgroup).

\begin{lemma} 
Assume Hypothesis \ref{auto} and let $(t,z)\in\bbR\times\bbC_+$. \\
$(i)$ If $\det(X)>0$, then
\begin{equation}\label{>0}
g_t(z)=
\frac{\big (\cos(\omega t)+\frac{\beta}{\omega} \sin
(\omega t)\big ) z+ \frac{\alpha}{\omega}\sin
(\omega t) }
{\big ( \frac{\gamma}{\omega}\sin
(\omega t) \big )z+\cos(\omega t)-\frac{\beta}{\omega} \sin
(\omega t)},
\end{equation}
where $\omega=\sqrt{\det(X)}>0$. \\
$(ii)$ If $\det(X)=0$, then
\begin{equation}\label{=0}
g_t(z)=
\frac{\big (1+\beta t\big ) z+ 
\alpha t }
{ \gamma t z+(1-\beta t)}.
\end{equation}
$(iii)$ If $\det(X)<0$, then
\begin{equation}\label{<0}
g_t(z)=
\frac{\big (\cosh(\omega t)+\frac{\beta}{\omega} \sinh
(\omega t)\big ) z+ \frac{\alpha}{\omega}\sinh
(\omega t) }
{\big ( \frac{\gamma}{\omega}\sinh
(\omega t) \big )z+\cosh(\omega t)-\frac{\beta}{\omega} \sinh
(\omega t)},
\end{equation}
where $\omega=\sqrt{|\det(X)|}>0$.
\end{lemma}
\begin{proof}
Since $\tr(X)=0 $, every entry $a_t$, $b_t$, $c_t$, and $d_t$ of the
matrix $e^{tX}$ in \eqref{entry} is a solution of the initial value
problem
\begin{equation}\label{dif}
\ddot y+\det(X)y=0,
\end{equation}
\begin{equation}\label{I}
y(0)=1, \quad \dot y(0)=\beta \quad \text{ for } y(t)=a_t,
\end{equation}
\begin{equation}\label{II}
y(0)=0, \quad \dot y(0)=\alpha \quad \text{ for } y(t)=b_t,
\end{equation}
\begin{equation}\label{III}
y(0)=0, \quad \dot y(0)=\gamma \quad \text{ for } y(t)=c_t,
\end{equation}
 \begin{equation}\label{IV}
y(0)=1, \;\;\, \dot y(0)=-\beta \;\, \text{ for } y(t)=d_t,
\end{equation}
where the dot $\cdot$ denotes $d/dt$. Solving the initial value
problems \eqref{dif}, \eqref{I}--\eqref{IV} proves
\eqref{>0}--\eqref{<0}.
\end{proof}

\begin{remark} \lb{r2.3} 
If $\gamma=0$ in \eqref{2.3} the subgroup $g_t$ is a group
of linear transformations of $\bbC_+$. If
$\gamma\in\bbR\backslash\{0\}$,
the subgroup $g_t$ corresponds to the case of linear fractional
transformations of $\bbC_+$.
If $\gamma\in\bbR\backslash\{0\}$, the automorphism $g_t(z)$ is a linear
function in
$z$  if and only if $t\in (\pi/\omega)\bbZ$ in case I and  
 $t=0$ in cases II and III, respectively. In other words, 
\begin{equation}
\text{$\gamma\in\bbR\backslash\{0\}$ if and only if } c_t\neq 0 \text{
for
}\begin{cases}
\text{$t\in \bbR\backslash\{(\pi/\omega)\bbZ\}$} & \text{in case I,} \\
\text{$t\in\bbR\backslash\{0\}$} & \text{in cases II, III.}
\end{cases} \lb{gam}
\end{equation}
Moreover, suppose that $t\in\bbR\backslash\{ (\pi/\omega)\bbZ\}$, that
is, $g_t$ is not the identity transformation $($$g_t(z)\neq z$$)$. Then
case I consists of elliptic M\"obius transformations. Case II always
corresponds to parabolic M\"obius transformations, and as long as
$t\neq 0$, case III corresponds to hyperbolic M\"obius
transformations. 
\end{remark}

\begin{remark} \lb{r2.4}
The case of self-adjoint rank-one perturbations $tP$ of self-adjoint
operators $A$ discussed in Lemma \ref{l1.1}, corresponds to the case
$\det(X)=0$ with $\alpha=\beta=0$, $\gamma=1$ as one readily verifies
upon comparison with \eqref{pert}. Similarly, the case of self-adjoint
extensions of a closed symmetric densely defined operator $\dot A$
with deficiency indices $(1,1)$ discussed in Lemma \ref{l1.3}
corresponds to the case $\det(X)=1$, $\omega=1$ with $\alpha=-1$,
$\beta=0$, $\gamma=1$ upon comparison with \eqref{ext} and the change
of parametrization $t\mapsto \tan(t)$ in \eqref{1.23}. 
\end{remark}

\begin{remark} \lb{r2.5} 
The geometry of the trajectories $\bigcup_{t\in \bbR}
\{g_t(z)\}$, $z\in\bbC_+$ of the one-parameter groups of automorphisms
\eqref{>0}--\eqref{<0} can be understood in terms of the trajectories
$\bigcup_{t\in \bbR} \{F_t(z)\}$, $z\in\bbC_+$, of the map $F_t$ given by
\begin{equation}
F_t(z)=\frac{\big (1+\beta t\big ) z+ 
\alpha t }
{ \gamma t z+(1-\beta t)}, \quad (\alpha, \beta, \gamma)\in \bbR^3.
\end{equation}
In fact, one has the following representations
$($$X=X(\alpha,\beta,\gamma)$, cf.\ \eqref{2.3}$)$
\begin{equation}
g_t(z)=\begin{cases}F_{\tan (\sqrt{\det(X)}t )/\sqrt{\det(X)}}(z)
&\text{ if }\det(X)=-\alpha\gamma-\beta^2>0,\\
F_t(z)&\text{ if }\det(X)=-\alpha\gamma-\beta^2=0,\\
F_{\tanh (\sqrt{|\det(X)|}t)/\sqrt{|\det(X)|}}(z)&\text{ if
}\det(X)=-\alpha\gamma-\beta^2<0. \end{cases}
\end{equation}
Therefore, the trajectories of the groups \eqref{>0}--\eqref{<0}
can be described by
\begin{equation}
\bigcup_{t\in \bbR} \{g_t(z)\}=\bigcup_{t\in \bbR} \{F_t(z)\}
\quad \text{ in cases I and II}
\end{equation}
and
\begin{equation}
\bigcup_{t\in \bbR} \{g_t(z)\}=\bigcup_{|t|<|\det(X)|^{-1/2}} \{F_t(z)\}
\subsetneqq \bigcup_{t\in \bbR} \{F_t(z)\} \quad \text{ in case III}. 
\end{equation}
One observes that $F_t$ is a one-parameter group of transformations of 
$\bbC_+$ with respect to $t$, that is, $F_{t+s}=F_t\circ F_s$ for all
$s,t\in \bbR$,  if and only if $\alpha\gamma+\beta^2=0$.
\end{remark}

Next, we denote by $\log(\cdot)$ the branch of the logarithm on the cut
plane $\Pi=\bbC\backslash [0,\infty)$ assuming
\begin{equation}\label{stanbr}
0<\arg (\log(z))<2\pi \text{ for } z\in\Pi,
\end{equation}
extending $\log(\cdot)$ to the upper rim, $\partial_+\Pi$, of $\Pi$ by
\begin{equation}
\lim_{\varepsilon\downarrow 0}\log(x+i\varepsilon)\in \bbR, \quad x>0 
\lb{rim}
\end{equation}
and hence
\begin{equation}
\Im (\log(x))=\pi, \quad x<0.
\end{equation}
Analytic continuation of the the branch $\log (\cdot)$ defined above
then
leads to the infinitely sheeted Riemann surface $\cR$ of the logarithm
with branch points of infinite order at zero and infinity. We denote
the resulting analytic function on $\cR$ by $\Ln(\cdot)$. For future
reference we also introduce the $n$th sheet, $\cS_n$, of $\cR$. We use
the convention $\cS_0=\Pi\cup\partial_+\Pi$. $\Ln\colon v\mapsto
w=\Ln(v)$ then maps the interior$,\Int(\cS_n)$, of each sheet $\cS_n$
biholomorphically onto the strip $2\pi n<\Im(z)<2\pi (n+1)$ and  
\begin{equation}
v\in\cS_n \text{ if and only if } 2\pi n \leq \arg(w) < 2\pi (n+1),
\quad
n\in\bbZ.
\end{equation}
Assuming Hypothesis \ref{auto} with $\gamma\neq 0$, we will in the
following 
\begin{align}
\begin{split}
& \text{denote by $\overleftrightarrow{c_tz+d_t}$ the lift of the
trajectory
$t\mapsto c_tz+d_t$ to $\cR$,} \\ 
& \text{with $\overleftrightarrow{c_0 z+d_0}=d_0=1\in \partial\cS_0$, 
\, $z\in\bbC_+$.}
\lb{lift}
\end{split}
\end{align}

\begin{lemma}\label{index}
Assume Hypothesis \ref{auto} with $\gamma\in\bbR\backslash\{0\}$ $($cf.\
\eqref{gam}$)$, let $-\infty<t_1<t_2 <\infty$, $z\in\bbC_+$, and recall
our convention \eqref{lift}. Then, 
\begin{equation}\label{IM}
\int_{t_1}^{t_2} \, dt \,\,\Im (g_t(z))=\frac{1}{\gamma}
\Im \big (\Ln(\overleftrightarrow{c_{t_2}
z+d_{t_2}})-\Ln(\overleftrightarrow{c_{t_1}z+d_{t_1}})\big ).
\end{equation}
\end{lemma}
\begin{proof} Since the entries of the matrix \eqref{entry} solve the
system of differential equations
\begin{equation}
\frac{d}{dt}\begin{pmatrix} a_t & b_t \\ c_t & d_t\end{pmatrix}=
\begin{pmatrix} \beta & \alpha \\ \gamma & -\beta \end{pmatrix}
\begin{pmatrix} a_t & b_t \\ c_t & d_t \end{pmatrix}, \quad t\in\bbR, 
\end{equation}
the following relations hold
\begin{equation}
\dot c_t=\gamma a_t-\beta c_t, \quad \dot d_t=\gamma b_t -\beta d_t, 
\end{equation}
implying
\begin{equation}
a_t=\frac{\dot c_t+\beta c_t}{\gamma} \quad \text{ and }\quad  
b_t=\frac{\dot d_t+\beta d_t}{\gamma}. 
\end{equation}
Thus, 
\begin{equation}
a_tz+b_t=\frac{\beta}{\gamma}(c_tz+d_t)
+\frac{1}{\gamma}(\dot c_tz+\dot d_t), 
\end{equation}
and hence
\begin{equation}\label{Ln}
g_t(z)=\frac{a_tz+b_t}{c_tz+d_t}=\frac{\beta}{\gamma}
+\frac{1}{\gamma}\frac{\dot c_tz+\dot
d_t}{c_tz+d_t}=\frac{\beta}{\gamma}+
\frac{1}{\gamma}\frac{d}{dt} \Ln(\overleftrightarrow{c_{t}z+d_{t}}).
\end{equation}
Integrating \eqref{Ln} from $t_1$ to $t_2$ and taking imaginary
parts of the resulting expression proves \eqref{IM}.
\end{proof}

\section{Dynamical systems on a space of measures} \lb{s3}

As shown below, each one-parameter subgroup $\{e^{tX}\}_{t\in \bbR}$ of
$SL_2(\bbR)$, or, what is the same, each one-parameter group
$\{g_t\}_{t\in \bbR}$ of automorphisms of the open upper-half plane
$\bbC_+$, generates  a dynamical system $\{g_t^*\}_{t\in \bbR}$ on the
(phase)  space
${\cM}=[0,\infty)\times\bbR\times \Omega$. Here $\Omega$ denotes
the  space of Borel measures $\mu$ on $\bbR$ with the property
\begin{equation}
\int_\bbR\frac{d\mu(\lambda)}{1+\lambda^2}<\infty. 
\end{equation}

Let $\{ g_t\}_{t\in \bbR}$ be a one-parameter subgroup of
$\Aut(\bbC_+)$,
the group of automorphisms  of $\bbC_+$, 
\begin{equation}\label{wagen}
g_t(z)=\frac{a_tz+b_t}{c_tz+d_t},  \quad (t,z)\in \bbR\times\bbC_+.
\end{equation}
Given a point $(A_0,B_0,\mu_0)\in {\cM}$, introduce the Herglotz
function 
\begin{equation}\label{init}
M_0(z)=A_0z+B_0+\int_\bbR d\, \mu_0(\lambda)\bigg
(\frac{1}{\lambda-z}-\frac{\lambda}{1+\lambda^2}\bigg ), \quad z\in
\bbC_+,
\end{equation}
where
\begin{equation}
\mu_0((\lambda_1,\lambda_2))+\f{1}{2}\mu_0(\{\lambda_1\})
+\f{1}{2}\mu_0(\{\lambda_2\})=\f{1}{\pi}\lim_{\varepsilon\downarrow 0}
\int_{\lambda_1}^{\lambda_2} d\lambda \, \Im(M_0(\lambda+i\varepsilon)). 
\lb{mu}
\end{equation}
Since for each $t\in \bbR$, $g_t\in\Aut(\bbC_+)$, the one-parameter
family of functions 
\begin{equation}\label{hhhh}
M_t(z)=g_t(M_0(z)), \quad (t,z)\in \bbR\times\bbC_+
\end{equation}
is a one-parameter family of Herglotz functions. Therefore, $M_t$
admits the representation
\begin{equation}\label{hrt}
M_t(z)=A_tz+B_t+\int_\bbR d\, \mu_t(\lambda)\bigg
(\frac{1}{\lambda-z}-\frac{\lambda}{1+\lambda^2}\bigg ), \quad 
(t,z)\in\bbR\times\bbC_+
\end{equation}
for a unique triple $(A_t, B_t, \mu_t)\in {\cM}$. Define the map
\begin{equation}
g^*_t:\cM\to \cM, \quad 
(A_0, B_0, \mu_0)\mapsto (A_t, B_t, \mu_t), \quad t\in
\bbR. 
\end{equation}
Then, 
\begin{equation}
g^*_{t+s}=g^*_{t}\circ g^*_{s},\quad s,t \in \bbR. 
\end{equation}
That is, $\{g_t^*\}_{t\in \bbR}$ defines a dynamical system  on 
${\cM}$ as claimed.

We note that 
\begin{equation}
M_t(i)=A_ti+B_t+i\int_\bbR  \frac{d\,\mu_t(\lambda)}{1+\lambda^2}
\end{equation}
and thus, 
\begin{equation}
A_t=\int_\bbR  \frac{d\,\mu_t(\lambda)}{1+\lambda^2}-\Im(M_t(i)), 
\quad B_t=\Re(M_t(i)),\quad t\in \bbR. 
\end{equation}
Moreover, if $c_t\ne 0$ in \eqref{wagen}, then $A_t=0$, and hence
\begin{equation}
 \int_\bbR  \frac{d\,\mu_t(\lambda)}{1+\lambda^2}=\Im(M_t(i))
\text{ if } c_t\ne 0. 
\end{equation}

For the remainder of this section it is convenient to introduce the
following assumptions.

\begin{hypothesis} \lb{h3.1} Assume Hypothesis \ref{auto} and
\begin{equation}
\text{$\gamma \in \bbR\backslash\{0\}$, or equivalently, } 
c_t\neq 0 \text{ for } \begin{cases}
\text{$t\in \bbR\backslash\{(\pi/\omega)\bbZ\}$} & \text{in case I,}
\\
\text{$t\in\bbR\backslash\{0\}$} & \text{in cases II, III.}
\end{cases} \lb{gamma}
\end{equation}
\end{hypothesis}

The following statement is a variant of the exponential
Herglotz representation theorem due to Aronszajn-Donoghue \cite{AD56}
(see
also \cite{AD64}).

\begin{lemma}
Assume Hypothesis \ref{h3.1}, let $(z,t)\in\bbC_+\times\bbR$, and
recall our convention \eqref{lift}. Given a Herglotz function $M_0$
with  $M_0(i)\ne 0$, introduce the function 
\begin{equation}
N_t(z)= \Ln(\overleftrightarrow{c_{t} M_0(z)+d_{t}}). 
\end{equation} 
Then $N_t(\cdot)$ is  analytic on $\bbC_+$ and the following
representation holds 
\begin{equation}\label{state}
N_t(z)=\Re (N_t(i))+\int_\bbR d\lambda \, \xi_t(\lambda)\bigg
(\frac{1}{\lambda-z}-\frac{\lambda}{1+\lambda^2}\bigg ), 
\end{equation}
where
\begin{equation}
\xi_t(\lambda)=\frac{1}{\pi}\lim_{\varepsilon \downarrow 0}
\Im\big (\Ln(\overleftrightarrow{c_{t}M_0(\lambda+i
\varepsilon)+d_t})\big ) 
 \text{ for a.e$.\lambda\in\bbR$.} \lb{XI}
\end{equation}
\end{lemma}
\begin{proof}
Since $M_0(i)\ne 0$, the expression $c_t M_0(z)+d_t$ never vanishes, and
hence the lift $\overleftrightarrow{c_{t}M_0(z)+d_t}$ is well-defined
as a point on $\cR$. To set the stage, we assume that
$\overleftrightarrow{c_{t}M_0(z)+d_t}$ is a point  on the $n$th sheet
$\cS_n$ of $\cR$ for some (and hence for all) $z\in\bbC_+$, that is, 
\begin{equation}
2\pi n\le \arg (\overleftrightarrow{c_{t}M_0(z)+d_{t}}) <2\pi (n+1),
\quad n\in
\bbZ. 
\end{equation}
Then, by the definition of $\Ln(\cdot)$ on $\cR$, one obtains
\begin{equation}
N_t(z)=\log(c_tM_0(z)+d_t)+2\pi i n, 
\end{equation}
where $\log (\cdot)$ denotes the branch \eqref{stanbr}, \eqref{rim} on
$\cS_0=\Pi\cup\partial_+\Pi$.

Given $t\in \bbR$, there are three possible outcomes for $N_t$ depending
on whether $c_t>0$, $c_t<0$, and $c_t=0$. If $c_t>0$, the function
$c_tM_0(z)+d_t$ is a Herglotz function and thus,
\begin{equation}\label{onsheet}
N_t(z)=\Re (N_t(i))+\int_\bbR d\lambda \, \eta_t(\lambda)\bigg
(\frac{1}{\lambda-z}-\frac{\lambda}{1+\lambda^2}\bigg )+2\pi i n,
\end{equation}
where
\begin{equation}
\eta_t(\lambda)=\frac{1}{\pi}\Im\big (\log( c_t m_0(\lambda)+d_t)\big )
\text{ for a.e$.\lambda\in\bbR$}
\end{equation}
and 
\begin{equation}
m_0(\lambda)=\lim_{\varepsilon \downarrow 0}M_0(\lambda+i\varepsilon)
\text{ for a.e$.\lambda\in\bbR$.} \lb{3.20}
\end{equation}
Since 
\begin{equation}
\frac{1}{\pi}\int_\bbR d\lambda \, \bigg
(\frac{1}{\lambda-z}-\frac{\lambda}{1+\lambda^2}\bigg )=i, \quad
z\in\bbC_+, 
\end{equation}
one can rewrite \eqref{onsheet} in the form \eqref{state}
with
\begin{align}
\xi_t(\lambda)&=\eta_t(\lambda)+2 n
\no \\
&=\frac{1}{\pi}
\Im\big (\log( c_t m_0(\lambda)+d_t)+2\pi n i \big )
\no \\
&=\frac{1}{\pi}
\lim_{\varepsilon \downarrow 0}\Im\big (\Ln(
\overleftrightarrow{c_{t}M_0(\lambda+i\varepsilon)+d_t})\big ),
\end{align}
proving \eqref{state}, \eqref{XI} in the case $c_t>0$. 

If $c_t<0$ one obtains 
\begin{align}
N_t(z)&=\log(c_tM_0(z)+d_t)+2\pi i n \no \\
&=\log(|c_t|M_0(z)-d_t)+2\pi i n +\pi i.
\end{align}
Using the Herglotz representation theorem for $|c_t|M_0(z)-d_t$ one
arrives at 
\begin{equation}
N_t(z)=\Re (N_t(i))+\int_\bbR d\lambda \, \eta_t(\lambda)\bigg
(\frac{1}{\lambda-z}-\frac{\lambda}{1+\lambda^2}\bigg )
+2\pi i n+\pi i, 
\end{equation} 
where
\begin{equation}
\eta_t(\lambda)=\frac{1}{\pi}
\Im (\log(  |c_t|m_0(\lambda)-d_t)) \text{ for a.e$.\lambda\in\bbR$}
\end{equation}
and \eqref{3.20} results again. Thus, \eqref{state} holds with
\begin{align}
\xi_t(\lambda)&=\eta_t(\lambda)+2 n -1=
\no \\
&=\frac{1}{\pi}\lim_{\varepsilon \downarrow 0}
\Im\big (\log(c_t M_0(\lambda+i\varepsilon)+d_t)+2\pi n i \big )
\no \\
&=\frac{1}{\pi}\Im\big
(\Ln(\overleftrightarrow{c_{t}m_0(\lambda)+d_t})\big ).
\end{align}
Finally, if $c_t=0$, $N_t(z)$ is a constant with respect to $z$
\begin{align}
N_t(z)&=\log (d_t)+2\pi i (n-1)
\no \\
& =\log |d_t|+i(2\pi (n-1)+\arg(d_t)),\quad \Im(z)>0,
\end{align}
which proves \eqref{state} with
\begin{equation}
\xi_t(\lambda)=2(n-1)+\pi^{-1}\arg(d_t)\in\bbZ, 
\end{equation}
a $\lambda$-independent integer constant.
\end{proof}

Now we can prove the absolute continuity of the measure associated with 
the Herglotz representation of the integrated (averaged) Herglotz
function
\begin{equation}
M_{t_1,t_2}(z)=\int_{t_1}^{t_2} dt \, M_t(z), \quad \Im(z)>0, \quad
t_1, t_2\in\bbR, \; t_1<t_2. 
\end{equation}

\begin{theorem}\label{mz}
Assume Hypothesis \ref{h3.1}, let $z\in\bbC_+$, $t_j\in\bbR$, $j=1,2$,
$t_1<t_2$, and recall our convention \eqref{lift}. Then the integrated
Herglotz function
\begin{equation}\label{integ}
M_{t_1,t_2}(z)=\int_{t_1}^{t_2} dt \, M_t(z)
\end{equation}
admits the  Herglotz representation
\begin{equation}\label{integm}
M_{t_1,t_2}(z)=B_{t_1,t_2}+\int_\bbR d\mu_{t_1,t_2}(\lambda)\bigg
(\frac{1}{\lambda-z}-\frac{\lambda}{1+\lambda^2} \bigg ). 
\end{equation}
Here $B_{t_1,t_2}\in \bbR$ and the measure $\mu_{t_1,t_2}$ is
absolutely continuous with respect to Lebesgue measure on $\bbR$ with
Radon--Nikodym derivative
$($density$)$ a bounded function
$\frac{d\mu_{t_1,t_2}}{d\lambda}=\xi_{t_1,t_2}\in L^\infty(\bbR)$. In
fact, 
\begin{equation}\label{xi}
\xi_{t_1,t_2}(\lambda)=\frac{1}{\gamma}\big
(\xi_{t_2}(\lambda)-\xi_{t_1} (\lambda)\big),
\end{equation}
where 
\begin{equation}
\xi_t(\lambda)=\lim_{\varepsilon \downarrow 0}\frac{1}{\pi}
\Im \big(\Ln 
(\overleftrightarrow{c_{t}M_0(\lambda+i\varepsilon)+d_t})\big), 
\quad t\in \bbR. 
\end{equation}
\end{theorem}
\begin{proof} By Lemma \ref{index}, $\Im(M_{t_1,t_2})$
admits the representation
\begin{align}
\Im (M_{t_1,t_2}(z))&=
\int_{t_1}^{t_2} \, dt \, \Im (g_t(M_t(z)))
\no \\
&=\frac{1}{\gamma}
\Im \big (\Ln(\overleftrightarrow{c_{t_2}M_0(z)+d_{t_2}})
-\Ln(\overleftrightarrow{c_{t_1}M_0(z)+d_{t_1}})\big ).
\label{hrepr}
\end{align}
Hence, $\Im(M_{t_1,t_2})$ is uniformly bounded on $\bbC_+$, which
proves that $M_{t_1,t_2}$ has no linear term in its Herglotz 
representation. Moreover, by Fatou's theorem, the boundedness of 
$\Im(M_{t_1,t_2})$ on $\bbC_+$ ensures the absolute continuity of the
measure $\mu_{t_1,t_2}$ in \eqref{integm} with respect to Lebesgue
measure on $\bbR$. Hence, \eqref{xi} is a consequence of \eqref{hrepr}.
\end{proof}

\begin{corollary}\label{corr}
Assume in addition to the hypotheses of Theorem \ref{mz} that
$\gamma>0$, and $t_1 <0<t_2$. If $\det(X)>0$, assume in addition that 
\begin{equation}
-\frac{\pi}{2 \sqrt{\det(X)}}<t_1 <0< t_2<\frac{\pi}{2 \sqrt{\det(X)}}. 
\end{equation}
Then the density \eqref{xi}
has the form
\begin{equation}
\xi_{t_1,t_2}(\lambda)=\frac{1}{\gamma} +\frac{1}{\gamma \pi}\Im \bigg
(\log \bigg(\frac{\Theta(t_2) m_0(\lambda)+1}{-\Theta(t_1)
m_0(\lambda)-1}\bigg)\bigg ) \text{ for a.e.\ $\lambda\in\bbR$,}
\end{equation}
where
\begin{equation}
m_0(\lambda)=\lim_{\varepsilon \downarrow 0} (\gamma
M_0(\lambda+i\varepsilon)-\beta) \text{ for a.e$.\lambda\in \bbR$}, 
\lb{3.37} 
\end{equation} 
and 
\begin{equation}
\Theta(t)=\lim_{s\to \sqrt{\det(X)}} \frac{\tan (st )}{s}
=\begin{cases}\frac{\tan (\sqrt{\det(X)}t )}{\sqrt{\det(X)}}
&\text{ if }\det(X)>0,\\
t&\text{ if }\det(X)=0,\\
\frac{\tanh (\sqrt{|\det(X)|}t)}{\sqrt{|\det(X)|}}&\text{ if }\det(X)<0,
\end{cases} \quad t\in\bbR. 
\end{equation}
\end{corollary} 

\begin{remark}\label{remcorr}
Define
\begin{equation}
T_1(X)=\begin{cases} -\frac{\pi}{2 \sqrt{\det(X)}}, & \det(X)>0, \\
-\infty, & \det(X) \le 0, \end{cases}, \quad 
T_2(X)=\begin{cases} \frac{\pi}{2 \sqrt{\det(X)}}, & \det(X)>0, \\
\infty, & \det(X) \le 0, \end{cases}, \lb{3.40} 
\end{equation}
then the density \eqref{xi} has the form
\begin{equation}
\xi_{T_1(X),T_2(X)}(\lambda)=\frac{1}{\gamma}+\begin{cases} 0 &\text{
if }\det(X)\ge 0, \\ 
\frac{1} {\gamma\pi}
\Im \bigg ( \log\bigg( \frac{\pi m_0(\lambda)+2\sqrt{|\det(X)|}}{\pi
m_0(\lambda)-2\sqrt{|\det(X)|}}\bigg )\bigg) &\text{ if }\det(X)<0.
\end{cases}
\end{equation}
\end{remark}

Next, we discuss the following technical result. 

\begin{lemma}\label{fubini}
Assume Hypothesis \ref{h3.1}, let $t_j\in\bbR\cup\{-\infty,\infty\}$,
$j=1,2$,
$t_1<t_2$, and denote by $\mu_{t_1,t_2}$ the Borel
measure in the Herglotz representation \eqref{integm} of the
integrated Herglotz function \eqref{integ}. Then for any bounded Borel
set $\Delta\subset
\bbR$ the function $t\mapsto \mu_t (\Delta)$ is measurable and one has 
\begin{equation}\label{aver1}
\int_{t_1}^{t_2} dt \, \mu_t (\Delta)=\mu_{t_1,t_2}(\Delta).
\end{equation}
\end{lemma}
\begin{proof}
The proof is based on the following representation
\begin{equation}\label{aver}
\int_{t_1}^{t_2} dt
\int_{\bbR} d\mu_t (\lambda) \, f(\lambda)=
\int_{\bbR}d\mu_{t_1,t_2}(\lambda) \, f(\lambda),
\end{equation}
which holds for a wide function class of $f$ to be specified below. 
We split the proof into four steps. First, we establish \eqref{aver}
for functions $f$ of the following type
\begin{equation}\label{step11}
f(\lambda)=\phi_{\varepsilon}(\lambda-\delta), \quad
\varepsilon >0, \quad \delta\in \bbR, 
\end{equation}
where $\phi_\varepsilon (\lambda)$ is an approximate identity, 
\begin{equation}\label{step12}
\phi_\varepsilon
(\lambda)=\varepsilon^{-1}\phi(\varepsilon^{-1}\lambda), \text{ with } 
\phi(\lambda)=\frac{1}{\pi}\frac{1}{1+\lambda^2}, \quad \lambda\in\bbR.
\end{equation}
Second, we prove 
\eqref{aver} for the functions $f$
that can be represented as a convolution of 
$\phi_\varepsilon  $ with a $C^\infty_0$-function $k$, 
\begin{equation}\label{step2}
f(\lambda)=(\phi_\varepsilon * k)(\lambda), \quad 
 k\in C_0^\infty(\bbR), \quad \varepsilon >0.
\end{equation}
Third, we prove the validity of representation \eqref{aver} for $f\in
C^\infty_0(\bbR)$. Finally,  we establish \eqref{aver} for
characteristic
functions of  finite intervals, implying assertion \eqref{aver1}.

{\it Step I.}
Let $z=\delta+i\varepsilon \in \bbC_+$. By representation \eqref{hrt}
\begin{equation}\label{hhh1}
\Im \big ( g_t(M_0(\delta+i\varepsilon))\big )=
A_t \varepsilon +\varepsilon \int_{\bbR}\frac{d\mu_t(\lambda)}
{(\lambda-\delta)^2+\varepsilon^2}.
\end{equation}
Since $\gamma\neq 0$, $A_t=0$ for almost all $t\in \bbR$ by \eqref{gam}
and hence
\begin{equation}\label{r1}
\int_{t_1}^{t_2} dt \, \Im \big (
g_t(M_0(\delta+i\varepsilon))\big )=\varepsilon \int_{t_1}^{t_2}
dt\int_{\bbR}\frac{d\mu_t(\lambda)} {(\lambda-\delta)^2+\varepsilon^2}.
\end{equation}
On the other hand, by Theorem \ref{mz} one infers
\begin{equation}\label{r2}
\int_{t_1}^{t_2} dt \, \Im \big (
g_t(M_0(\delta+i\varepsilon))\big) =\varepsilon
\int_{\bbR}\frac{d\mu_{t_1,t_2}(\lambda)
}{(\lambda-\delta)^2+\varepsilon^2}.
\end{equation}
Comparing \eqref{r1}  and \eqref{r2} proves \eqref{aver}
for the functions of the type \eqref{step11},
\eqref{step12}.

{\it Step II}. Let $ k\in C_0^\infty(\bbR)$ and 
$\supp ( k)\subset [\delta_1, \delta_2]$ 
for some $-\infty < \delta_1 < \delta_2 <\infty$.

We start with two observations. Given $\varepsilon >0$, the function
$k_\varepsilon(\lambda, \delta)=\phi_{\varepsilon}(\lambda-\delta)
k(\delta)$, $(\lambda,\delta)\in\bbR\times [\delta_1,\delta_2]$, is
summable with respect to the product measures $d\mu_t\times
d\delta$, $t\in [t_1, t_2]$ as well as with respect to the
product measure $d\mu\times d\delta$, that is,
\begin{equation}\label{sum1}
k_\varepsilon \in L^1(\bbR\times[\delta_1,
\delta_2]; d\mu_t\times d \delta),\quad t\in [t_1, t_2]
\end{equation}
and
\begin{equation}\label{sum2a}
k_\varepsilon \in L^1(\bbR\times[\delta_1,\delta_2]; d\mu\times d
\delta).
\end{equation}
Moreover, we claim that the function $F_\varepsilon$ (cf.\
\eqref{step2}),
\begin{equation}\label{sum}
F_\varepsilon(t, \delta)=\int_{\bbR}d\mu_t(\lambda)\,\phi_{\varepsilon}
(\lambda-\delta)f(\delta), \quad t\in [t_1,t_2], 
\end{equation}
is summable on  $[t_1, t_2]\times [\delta_1, \delta_2]$, that is,
\begin{equation}\label{sum2}
F_\varepsilon\in L^1([t_1, t_2]\times [\delta_1,
\delta_2];dt\times d\delta).
\end{equation}
In order to prove \eqref{sum2}, one notes that by representation
\eqref{hrt}, \eqref{hhh1} holds again. Thus, given $\varepsilon >0$,
the function $h$, 
\begin{equation}
h(t, \delta)=A_t \varepsilon +\varepsilon
\int_{\bbR}\frac{d\mu_t(\lambda)} {(\lambda-\delta)^2+\varepsilon^2}, 
\end{equation}
is continuous on $[t_1, t_2]\times [\delta_1, \delta_2]$. Hence 
$h(t,\delta)f(\delta)$ is also continuous on 
$[t_1, t_2]\times [\delta_1, \delta_2]$ and thus bounded.
Since $A_t=0$ a.e., $F_\varepsilon(t, \delta)$ is measurable and
essentially bounded on $[t_1, t_2]\times [\delta_1,\delta_2]$.
This proves \eqref{sum2}.

Now the validity of \eqref{aver} for the function class \eqref{step2}
follows from the following chain of equalities
\begin{align}
&\int_{t_1}^{t_2} 
dt\int_\bbR d\mu_t\,(\lambda)\,(\phi_\varepsilon * k )(\lambda) \no \\
&=\int_{t_1}^{t_2} dt\int_\bbR
d\mu_t(\lambda)\int_{\delta_1}^{\delta_2} d\delta \,
\phi_\varepsilon(\lambda-\delta) k (\delta)
\quad (\text{since } \supp(k)\subset [\delta_1, \delta_2]) \no \\
&=\int_{t_1}^{t_2} dt\int_{\delta_1}^{\delta_2}d\delta 
\int_\bbR d\mu_t(\lambda) \, \phi_\varepsilon(\lambda-\delta) k (\delta)
\quad (\text{by  \eqref{sum1} using Fubini's theorem)} \no \\
&=\int_{\delta_1}^{\delta_2} d\delta\int_{t_1}^{t_2} dt
\int_\bbR d\mu_t(\lambda) \, \phi_\varepsilon(\lambda-\delta) k (\delta)
\quad (\text{by  \eqref{sum2} using Fubini's theorem)} \no \\
&=\int_{\delta_1}^{\delta_2} d\delta\int_\bbR 
d\mu_{t_1,t_2}(\lambda)\,\phi_\varepsilon(\lambda-\delta) k (\delta)
\quad (\text{by  step I)} \no \\
&=\int_\bbR d\mu_{t_1,t_2}(\lambda) \int_{\delta_1}^{\delta_2} d\delta
\, \phi_\varepsilon(\lambda-\delta) k (\delta)
\quad (\text{by \eqref{sum2a} using Fubini's theorem)} \no \\
&=\int_\bbR d\mu_{t_1,t_2}(\lambda)\,(\phi_\varepsilon * k )(\lambda)
\quad (\text{since } \supp(k)\subset [\delta_1, \delta_2]).
\end{align}

{\it Step III.} Let $f\in C_0^\infty(\bbR)$ with $\supp(f)\subset
[\delta_1, \delta_2]$. One infers
\begin{equation}\label{unif} 
\lim_{\varepsilon \downarrow 0}(\phi_\varepsilon * f)(\lambda)=
\lim_{\varepsilon \downarrow 0}\int_{\delta_1}^{\delta_2}d \delta \, 
\phi_{\varepsilon}(\lambda-\delta)f(\delta)=f(\lambda), 
\end{equation}
uniformly with respect to $\lambda$ as long as $\lambda$ varies in a 
compact set $\Lambda\subset\bbR$. With 
$\Lambda=[\delta_1-1,\delta_2+1]$, one obtains the estimate
\begin{equation}
\bigg \vert \int_{\delta_1}^{\delta_2}d \delta \, 
\phi_{\varepsilon}(\lambda-\delta)f(\delta)\bigg \vert\le
\max_{\delta\in \supp(f)} |f(\delta)|(\delta_2-\delta_1)\frac{1}{\pi}
\frac{\varepsilon} {\text{ dist}^2(\lambda,[\delta_1, \delta_2])}, 
\quad \lambda\in \bbR\backslash \Lambda. 
\end{equation}
Thus, there exists a constant $C=C(\delta_1, \delta_2)$ such that
\begin{equation}\label{tails}
|(\phi_\varepsilon * f)(\lambda)|=
\bigg \vert \int_{\delta_1}^{\delta_2}d \delta \, 
\phi_{\varepsilon}(\lambda-\delta)f(\delta)\bigg \vert\le
C\frac{\varepsilon}{1+\lambda^2}, \quad \lambda\in \bbR\backslash
\Lambda.
\end{equation}
Taking into account that
\begin{equation}
\sup_{t\in
[t_1,t_2]}\int_\bbR\frac{d\mu_t(\lambda)}{1+\lambda^2}<\infty, 
\end{equation}
the uniform convergence \eqref{unif} combined with the estimate 
\eqref{tails} and the result of Step II proves \eqref{aver} for $f\in
C_0^\infty(\bbR)$.

{\it Step IV.}
Let $\Delta$ be a finite interval and $f_1(\lambda)\ge f_2(\lambda) \ge
\dots $ a monotone sequence of non-negative functions, $f_n\in
C_0^\infty(\bbR)$, $ n\in\bbN$ converging pointwise to the
characteristic function of the interval $\Delta$ as $n$  approaches
infinity, that is 
\begin{equation}
\lim_{n\to \infty } f_{n}(\lambda)=\chi_\Delta(\lambda),
\quad \lambda \in \bbR. 
\end{equation}
By the dominated convergence theorem one then obtains 
\begin{equation}\label{mud}
\lim_{n\to \infty}\int_{\bbR} d\mu_{t_1,t_2}(\lambda) \,
f_n(\lambda)=\int_{\bbR}d
\mu_{t_1,t_2}(\lambda)\,\chi_\Delta(\lambda)=\mu_{t_1,t_2}(\Delta)
\end{equation}
and
\begin{equation}
\lim_{n\to \infty}\int_{\bbR} d\mu_t(\lambda) \, f_n(\lambda)
=\int_{\bbR}d \mu_t(\lambda)\,\chi_\Delta(\lambda)=\mu_t(\Delta), \quad
t\in [t_1,t_2]. 
\end{equation}
Since
\begin{align}
0&\le\int_{\bbR} d\mu_t(\lambda) \, f_n(\lambda)\le
\int_{\bbR} d\mu_t(\lambda) \, f_1(\lambda) \no \\
& \le \max_{s\in \supp(f_1)} \big((1+s^2)f_1(s)\big) \int_{\bbR}
\frac{d\mu_t(\lambda)}{1+\lambda^2} \no \\
& \le \max_{s\in \supp(f_1)} \big((1+s^2)f_1(s)\big) \sup_{t\in
\bbR}\bigg(\int_{\bbR}\frac{d\mu_t(\lambda)}{1+\lambda^2}\bigg), 
\end{align}
one obtains
\begin{equation}
\lim_{n\to \infty}\int_{t_1}^{t_2}
dt \int_\bbR d\mu_t(\lambda) \, f_n(\lambda)=\int_{t_1}^{t_2}
dt \lim_{n\to \infty}\int_\bbR d\mu_t(\lambda) \, f_n(\lambda) 
=\int_{t_1}^{t_2}dt \, \mu_t(\Delta), 
\end{equation}
using the dominated convergence theorem again. By Step III and by 
taking into account \eqref{mud}, this proves \eqref{aver} for 
$f(\lambda)=\chi_\Delta(\lambda)$. 

The extension from the case of bounded intervals $\Delta$ to the case of
bounded Borel sets $\Delta$ is now straightforward, completing the
proof.
\end{proof}

Given a general Herglotz function $M$ of the type 
\begin{equation}
M(z)=Az+B+\int_\bbR d\mu(\lambda)\,\bigg(\f{1}{\lambda-z}
-\f{\lambda}{1+\lambda^2}\bigg), \quad A\geq 0, \; B\in\bbR,
\; z\in\bbC_+, \lb{hrhr} 
\end{equation}
we next introduce the
following  subsets  of $\bbR$,
\begin{align}
\cA(M)&=\Big \{\lambda\in \bbR\,\,\Big | \,\, \lim_{\varepsilon
\downarrow 0} M(\lambda+i\varepsilon)\in \bbC_+\Big \},
\label{invac} \\
\cP(M)&
=\Big \{\lambda\in \bbR \, \Big |\, \lim_{\varepsilon\downarrow
0}\varepsilon\,  
\Im( M(\lambda+i\varepsilon))\in (0,
\infty )\Big \}
\label{invpp}\\
&\quad \; \bigcup \Big \{\lambda\in \bbR \, \Big |\, \lim_{\varepsilon
\downarrow 0}\varepsilon^{-1} \Im( M(\lambda+i\varepsilon))\in (0,
\infty )\Big \},
\no \\
\cS(M)&=\bbR\backslash \{\cA(M)\cup\cP(M)\}.
\label{invsin}
\end{align}
By results of Aronszajn \cite{Ar57}, Donoghue \cite{Do65}, and Simon and
Wolff \cite{SW86}, the subsets
\eqref{invac}--\eqref{invsin} are 
 invariant with respect to the whole family $\{g\circ M\}_{g\in
\Aut(\bbC_+)}$
of Herglotz functions, that is, 
\begin{equation}
 \cA(M)=\cA(g\circ M),\quad \cP(M)=\cP(g\circ M),\quad \cS(M)=\cS(g\circ
M), \quad g\in \Aut(\bbC_+).
\end{equation}
Strictly speaking, these results were obtained for Herglotz functions
being the Stieltjes transforms of finite Borel measures. For the sake of
completeness we prove this invariance in the case of general Herglotz
functions. The invariance of the set $\cA(M)$ is obvious from
\eqref{>0}--\eqref{<0}. The invariance of the set $\cS(M)$ is then a
corollary of the one of
$\cP(M)$. In order to prove the invariance of the set $\cP(M)$ one needs
some additional considerations.

We start by recalling the following well-known result.

\begin{lemma} [see, e.g., \cite{Ar57}, \cite{AD56},
\cite{SW86}] \label{volk}
Let $M$ be a Herglotz function with representation \eqref{hrhr}. 
Then, for any $\lambda_0\in \bbR$, 
\begin{equation}
\lim_{\varepsilon\downarrow 0} (-i\varepsilon) 
M(\lambda_0+i\varepsilon)=\mu(\{\lambda_0\}), \lb{3.41A}
\end{equation}
in particular,
\begin{equation}
\lim_{\varepsilon\downarrow 0} \varepsilon 
\Im (M(\lambda_0+i\varepsilon))=\mu(\{\lambda_0\}) \lb{3.42}
\end{equation}
and
\begin{equation}
\lim_{\varepsilon\downarrow 0} \varepsilon 
\Re(M(\lambda_0+i\varepsilon))=0. \lb{3.43}
\end{equation}
\end{lemma}

\begin{definition}
A Herglotz function $M$ of the type \eqref{hrhr} is said to have a
normal derivative at the point $\lambda\in \bbR$ if the following two
limits exist $($finitely$)$. \\
$(i)$  $M(\lambda)=\lim_{\varepsilon \downarrow 0}
M(\lambda+i\varepsilon)\in\bbC$. \\
$(ii)$ $M'(\lambda)=\lim_{\varepsilon \downarrow 0} 
(M(\lambda+i\varepsilon)-M(\lambda))/(i\varepsilon)\in\bbC$. 
\end{definition}

\begin{lemma}\label{lisu}
Let $M$ be a Herglotz function with representation \eqref{hrhr}. 
Assume, in addition, that 
\begin{equation}\label{eq12}
\lim_{\varepsilon \downarrow 0}\varepsilon^{-1}\Im(
M(\lambda_0+i\varepsilon))\in (0, \infty)
\end{equation}
for some point $\lambda_0\in \bbR$. Then $M$ has a real normal boundary
value at $\lambda_0$ and $M$ has a strictly positive normal derivative
at $\lambda_0$, that is, 
\begin{align}
M(\lambda_0)&=\lim_{\varepsilon \downarrow 0}
M(\lambda_0+i\varepsilon)\in \bbR, \lb{3.73} \\
M'(\lambda_0)&=\lim_{\varepsilon \downarrow 0} 
\frac{M(\lambda_0+i\varepsilon)-M(\lambda_0)}{i\varepsilon}\in (0,\infty).
\label{eq11}
\end{align}
\end{lemma}
\begin{proof}
Let $\cI$ be a finite open interval containing $\lambda_0$ and decompose
$M$ as $M=M_1+M_2$, where
\begin{align}
M_1(z)&=Az+B+\int_{\bbR\backslash \cI}d\mu(\lambda)\bigg
(\frac{1}{\lambda-z}
-\frac{\lambda}{1+\lambda^2}\bigg )-\int_{\cI}d
\mu(\lambda)\frac{\lambda}{1+\lambda^2}, \label{inm1} \\
M_2(z)&=\int_{\cI}\frac{d\mu(\lambda)}{\lambda-z}. \label{inm2}
\end{align}
Clearly, 
\begin{equation}
M_1(\lambda_0)=\lim_{\varepsilon \downarrow 0}
M_1(\lambda_0+i\varepsilon)\in \bbR \lb{3.77}
\end{equation}
and
\begin{equation}\label{eq13}
M'_1(\lambda_0)=\begin{cases} \lim_{\varepsilon \downarrow 0} 
\frac{M_1(\lambda_0+i\varepsilon)-M_1(\lambda_0)}{i\varepsilon}>0  
& \text{if $A\neq 0$ or $\mu(\bbR\backslash \cI)\neq 0$,} \\
0 & \text{if $A=0$ and $\mu(\bbR\backslash \cI)=0$.} \end{cases}
\end{equation}
Hypothesis \eqref{eq12} and \eqref{inm2} imply
\begin{equation}
\lim_{\varepsilon \downarrow 0}\varepsilon^{-1}\Im(
M_1(\lambda_0+i\varepsilon))=0
\end{equation}
and hence
\begin{align}
\lim_{\varepsilon \downarrow 0}\varepsilon^{-1}\Im(
M(\lambda_0+i\varepsilon))&=
\lim_{\varepsilon \downarrow 0}\varepsilon^{-1}\Im(
M_2(\lambda_0+i\varepsilon))=\lim_{\varepsilon \downarrow 0}
\int_{\cI}\frac{d\mu(\lambda)}{(\lambda-\lambda_0)^2+\varepsilon^2} \no\\
&=\int_{\cI}\frac{d\mu(\lambda)}{(\lambda-\lambda_0)^2}\in [0, \infty), 
\label{eq14}
\end{align}
using the monotone convergence theorem in the last step. Since $\cI$ is a
finite interval and 
$\int_{\cI} d\mu(\lambda)(\lambda-\lambda_0)^{-2}<\infty$ by
\eqref{eq12} and \eqref{eq14}, 
applying the dominated convergence theorem yields
\begin{equation}
\lim_{\varepsilon \downarrow 0} 
\Re (M_2(\lambda_0+i\varepsilon))= \lim_{\varepsilon \downarrow 0}
\int_{\cI} d\mu(\lambda)
\frac{(\lambda-\lambda_0)} {(\lambda-\lambda_0)^2+\varepsilon^2}=
\int_{\cI}\frac{d\mu(\lambda)}{\lambda-\lambda_0} \in\bbR.  
\end{equation}
Thus,
\begin{equation}
M_2(\lambda_0)=\lim_{\varepsilon \downarrow 0}
M_2(\lambda_0+i\varepsilon)\in\bbR, \lb{3.81}
\end{equation}
and combining \eqref{3.77} and \eqref{3.81} then proves \eqref{3.73}.
Applying the dominated convergence theorem again yields
\begin{align}
M'_2(\lambda_0)&=\lim_{\varepsilon \downarrow 0}
\frac{M_2(\lambda_0+i\varepsilon)-M_2(\lambda_0)}{i\varepsilon}
= \lim_{\varepsilon \downarrow 0}\int_{\cI}
\frac{d\mu(\lambda)}{(\lambda-\lambda_0-i\varepsilon)(\lambda-\lambda_0)}
\no \\ 
&=\begin{cases} 
\int_{\cI}\frac{d\mu(\lambda)}{(\lambda-\lambda_0)^2}>0 & \text{if
$\mu(\cI)\neq 0$,} \\
0 & \text{if
$\mu(\cI) = 0$.} \end{cases}
\label{eq15}
\end{align}
Taking into account that by hypothesis \eqref{eq12}, either $A>0$ or
$\mu(\bbR)\ne 0$ in the Herglotz representation \eqref{hrhr} of $M$,
combining \eqref{eq13} and \eqref{eq15} proves \eqref{eq11}. 
\end{proof}

\begin{lemma}
Let $M$ be a Herglotz function of the type \eqref{hrhr}. Then
\begin{equation}
\cP(M)=\cP(g\circ M), \quad g\in \Aut(\bbC_+). \lb{3.57}
\end{equation}
\end{lemma}
\begin{proof}
It sufficies to prove the inclusion
\begin{equation}
\cP(M)\subset\cP(g\circ M), \quad g\in \Aut(\bbC_+). \lb{3.58}
\end{equation}
Moreover, any automorphism $g\in \Aut(\bbC_+)$ admits the representation
\begin{equation}
g=g_1\circ f \circ g_2, 
\end{equation}
where $g_j\in \Aut(\bbC_+)$, $j=1,2$ are linear transformations of the
upper half-plane, and
\begin{equation}
f(z)=-\frac{1}{z}, \quad z\in \bbC_+. 
\end{equation}
Since $\cP(M)$ is obviously invariant for linear transformations of
$\bbC_+$, it suffices to establish the inclusion
\begin{equation}
\cP(M)\subset \cP(f\circ M). \lb{3.59}
\end{equation} 
Let $\lambda\in \cP(M)$. By definition of $\cP(M)$ either
\begin{equation}\label{polus}
\lim_{\varepsilon \downarrow 0} 
\varepsilon\Im(M(\lambda+i\varepsilon))\in (0, \infty)
\end{equation}
or
\begin{equation}\label{nulili}
\lim_{\varepsilon \downarrow 0} 
\varepsilon^{-1} \Im(M(\lambda+i\varepsilon))\in (0, \infty).
\end{equation}
If \eqref{polus} holds, then
\begin{equation}
\lim_{\varepsilon\downarrow 0} \varepsilon^{-1} \Im( f\circ
M(\lambda+i\varepsilon))=
\lim_{\varepsilon\downarrow 0} \varepsilon^{-1}
\frac{\Im( M(\lambda+i\varepsilon))}{|M(\lambda+i\varepsilon)|^2}\in (0,
\infty),
\end{equation}
using \eqref{3.41A}--\eqref{3.43}. Therefore, $\lambda\in \cP(f\circ
M)$. Next, assume that \eqref{nulili} holds.
By Lemma \ref{lisu} $M(\lambda)=\lim_{\varepsilon \downarrow
0}M(\lambda+i\varepsilon)\in \bbR $  and $M(z)$ has a positive
normal derivative at the point $\lambda$. If $M(\lambda)\ne 0$ , then
\begin{equation}
\lim_{\varepsilon \downarrow
0}\varepsilon^{-1} \Im ((f\circ M)(\lambda+i\varepsilon))=
\lim_{\varepsilon \downarrow
0}\varepsilon^{-1}\Im \bigg(
\frac{1}{M(\lambda)}-\frac{1}{M(\lambda+i\varepsilon)}\bigg
)=\frac{M'(\lambda)}{(M(\lambda))^2}>0. 
\end{equation}
If $M(\lambda)=0$, then
\begin{equation}
\lim_{\varepsilon \downarrow
0}\varepsilon \Im ((f\circ M)(\lambda+i\varepsilon))=
\lim_{\varepsilon \downarrow
0}\varepsilon \Im\bigg (
-\frac{1}{M(\lambda+i\varepsilon)}\bigg
)=\frac{1}{M'(\lambda)}>0,   
\end{equation}
that is, \eqref{nulili} implies $\lambda\in\cP(f\circ M)$. Therefore, 
in both cases \eqref{polus} and \eqref{nulili}, $\lambda\in \cP(f\circ
M)$, which proves
\eqref{3.59} and hence \eqref{3.57}.
\end{proof}
 
The following result provides a spectral characterization of the
invariant sets $\cA(M)$, $\cS(M)$, and $\cP(M)$ (see \cite{SW86} for a
strategy of the proof). We recall that a measure $\mu$ on $\bbR$ is
supported by the set $\cT\subseteq\bbR$ if $\mu(\bbR\backslash\cT)=0$.
\begin{lemma}
Let $M$ be a Herglotz function of the type \eqref{hrhr}, $g\in
\Aut(\bbC_+)$, $\mu_g$ the measure in the Herglotz representation of
$g\circ M$, and
\begin{equation}
\mu_g=\mu_g^{ac}+\mu_g^{sc}+\mu_g^{pp}, \quad g\in \Aut(\bbC_+)
\end{equation}
the Lebesgue decomposition of $\mu_g$ into its absolute continuous, 
singularly continuous, and pure point parts, respectively.
Then $\mu_g^{ac}$, $\mu_g^{sc}$, and $\mu_g^{pp}$ are supported by
$\cA(M)$, $\cS(M)$ and $\cP(M)$, respectively. \\
Moreover, for any point $\lambda \in \cP(M)$ there exists an
automorphism $g\in \Aut(\bbC_+)$ such that 
\begin{equation}
\mu^{pp}_g(\{\lambda \})>0. 
\end{equation}
\end{lemma}

\begin{remark}
Originally, the set $\cP(M)$ has been introduced in the context of rank
one perturbations in \cite{SW86} by
\begin{align}
\cP(M)=
&\Big \{\lambda\in \bbR \, \Big|\, \mu(\{\lambda\})>0 \Big \}\bigcup
\Big \{\lambda\in \bbR \, \Big 
|\,\int_{\bbR}\frac{d\mu(s)}{(s-\lambda)^2}<\infty \Big \}.
\end{align}
\end{remark}

Naively one might think that the set $\bbR\backslash \cA(M)$ coincides
(modulo Lebesgue null sets)  with the complement of the support of the
absolutely continuous component $\mu^{\text{ac}}$ of the measure
$\mu$ associated with the Herglotz function $M$. Thus, one might
erroneously conclude that
\begin{equation}\label{mistake}
|\supp(\mu^{\text{ac}})\cap (\bbR\backslash \cA(M))|=0.
\end{equation}
The following counterexample illustrates the situation.

\begin{example}\label{exam}
Let  $ K\subset[0,1]$ be a closed nowhere dense set 
of a positive Lebesgue measure and let
\begin{equation}
M(z)=\int_{\bbR}\frac{d\mu(\lambda)}{\lambda-z}, \quad d\mu(\lambda)=
\chi_{[0,1]\backslash K}(\lambda)\,d\lambda, 
\end{equation}
where $\chi_\Lambda$ denotes the characteristic set of a set
$\Lambda\subset\bbR$. Then 
\begin{equation}
\supp(\mu)=
\supp(\mu^{\text{ac}})=[0,1], 
\end{equation}
but
\begin{equation}
|\supp(\mu^{\text{ac}})\cap (\bbR\backslash \cA(M))|=|K|>0. 
\end{equation}  
Thus, \eqref{mistake} does not hold in general.
\end{example}

Combining the results of Corollary \ref{corr}, Remark \ref{remcorr}, 
and Lemma \ref{fubini}, we can now formulate the following spectral
averaging theorems. 

\begin{theorem} \lb{t3.6}
Assume Hypothesis \ref{h3.1} and let
$t_j\in\bbR\cup\{-\infty,\infty\}$, $j=1,2$, $t_1<t_2$. Suppose $M_0$
is a Herglotz function of the type \eqref{hrhr} and $M_t$, $t\in\bbR$,
is the one-parameter family of Herglotz functions given by \eqref{hhhh}
and
\eqref{hrt}. Denote by $\mu^{\text{ac}}_t$, $\mu^{\text{sc}}_t$, and
$\mu^{\text{pp}}_t$ the absolutely continuous, singularly continuous,
and pure point parts in the Lebesgue decomposition of $\mu_t$ in
\eqref{hrt}, 
\begin{equation}\label{leb}
\mu_t=\mu^{\text{ac}}_t+\mu^{\text{sc}}_t+\mu^{\text{pp}}_t,\quad t\in
\bbR.
\end{equation}
Then the following averaged measures 
\begin{equation}
\int_{t_1}^{t_2} dt\, d\mu_t, \quad \int_{t_1}^{t_2}
dt\, d\mu_t^{\text{ac}}, \quad \int_{t_1}^{t_2}
dt\, d\mu_t^{\text{sc}}, \quad \int_{t_1}^{t_2}
dt\, d\mu_t^{\text{pp}} 
\end{equation}
are absolutely continuous with
respect to Lebesgue measure on $\bbR$. More precisely, given a  bounded
Borel set $\Delta\subset\bbR$, the functions $t\mapsto \mu_t(\Delta)$,
$t\mapsto \mu_t^{\text{ac}}(\Delta)$, 
$t\mapsto \mu_t^{\text{sc}}(\Delta)$, and 
$t\mapsto \mu_t^{\text{pp}}(\Delta)$  are measurable and
\begin{equation}
\int_{t_1}^{t_2} dt \begin{cases} \mu_t (\Delta) \\
\mu_t^{ac} (\Delta) \\ \mu_t^{sc} (\Delta) \\ \mu_t^{pp} (\Delta) 
\end{cases} \hspace{-3mm}= \,
\begin{cases}
\mu_{t_1,t_2}(\Delta) \\
\mu_{t_1,t_2}(\Delta\cap \cA) \\
\mu_{t_1,t_2}(\Delta\cap \cS) \\
\mu_{t_1,t_2}(\Delta\cap\cP) \\
\end{cases} \lb{3.61} 
\end{equation}
where $\mu_{t_1,t_2}$ is the absolutely continuous measure in the
Herglotz representation \eqref{integm} of the integrated Herglotz
function \eqref{integ} in Theorem \ref{mz} and $\cA(M_0)$, $\cP(M_0)$,
and $\cS(M_0)$ are the invariant  sets \eqref{invac}--\eqref{invsin}
associated with the Herglotz function $M_0$. In particular, 
\begin{align}
&|\{t\in \bbR\,\, |\,\, \mu_t^{\text{sc}}(\bbR)\ne 0\}|=0\quad 
\text{ if } |\cS(M_0)|=0,  \\
&|\{t\in \bbR\,\, |\,\, \mu_t^{\text{pp}}(\bbR)\ne 0\}|=0\quad 
\text{ if } |\cP(M_0)|=0.  
\end{align}
\end{theorem}
\begin{proof}
Equation \eqref{aver1} implies the result \eqref{3.61} since
\begin{align}
\mu_t(\Delta\cap \cA(M_0))=&\mu_t^{ac}(\Delta\cap \cA(M_0)),\no\\
\mu_t(\Delta\cap \cS(M_0))=&\mu_t^{sc}(\Delta\cap \cS(M_0)), \\
\mu_t(\Delta\cap\cP(M_0))=&\mu_t^{pp}(\Delta\cap\cP(M_0)), \no
\end{align}
and $\cA(M_0)$, $\cS(M_0)$, and $\cP(M_0)$ are known to be Borel sets.
\end{proof}

\begin{remark}
The ``life-time'' $|\{t\in \bbR \,\,|\,\, \mu_t^{\text{sing}}
(\bbR\backslash \cA(M_0))\ne 0\}|$ is never zero whenever $|\cS(M_0)\cup
\cP(M_0)|\ne 0$. Here
\begin{equation}
\mu_t^{sing}=\mu_t^{sc}+\mu_t^{pp}. 
\end{equation}
As concrete examples show $($cf.\ \cite{DFP01}$)$, it may be finite or
infinite depending upon the choice of the Herglotz function $M_0$. 
\end{remark}
\begin{remark}
Example \ref{exam} shows that the sets $\supp(\mu_t^{\text{ac}})$ and 
$\bbR\backslash \cA$ may have nontrivial intersection of positive
Lebesgue measure and that 
\begin{equation}
|\{t\in \bbR \,\,|\,\, \mu_t^{\text{sing}}
 (\supp(\mu_t^{\text{ac}}))\ne 0\}|\ne 0
\end{equation}
in general.
\end{remark}
As a corollary of the previous theorem
we get the following global result.

\begin{theorem}\label{globsa} 
Assume the hypotheses of Theorem \ref{t3.6} and let 
$\cA(M_0)$, $\cS(M_0)$, and $\cP(M_0)$ be the  invariant sets associated
with the Herglotz function $M_0$. Then for any bounded Borel set
$\Delta\subset
\bbR$ the following holds.\\
$(i)$ If $\det(X)>0$ then,
\begin{equation}
|\gamma|\int_{-\pi{\big/}\big(2\sqrt{\det(X)}\big)}^{\pi{\big/}
\big(2\sqrt{\det(X)}\big)} dt  \begin{cases} \mu_t (\Delta) \\
\mu_t^{ac} (\Delta) \\ \mu_t^{sc} (\Delta) \\ \mu_t^{pp} (\Delta) 
\end{cases} \hspace{-3mm}= \,
\begin{cases}
|\Delta| \\
|\Delta\cap \cA(M_0)| \\
|\Delta\cap \cS(M_0)| \\
|\Delta\cap\cP(M_0)| \\ \end{cases} \lb{3.107}
\end{equation}
$(ii)$ If $\det(X)=0$ then,
\begin{equation}
|\gamma|\int_{-\infty}^\infty dt  \begin{cases} \mu_t (\Delta) \\
\mu_t^{ac} (\Delta) \\ \mu_t^{sc} (\Delta) \\ \mu_t^{pp} (\Delta) 
\end{cases} \hspace{-3mm}= \,
\begin{cases}
|\Delta| \\
|\Delta\cap \cA(M_0)| \\
|\Delta\cap \cS(M_0)| \\
|\Delta\cap\cP(M_0)| \\ \end{cases} \lb{3.108}
\end{equation}
$(iii)$ If $\det(X)<0$ then,
\begin{align}
|\gamma|\int_{-\infty}^\infty dt \, \mu_t(\Delta)
&=\int_\Delta d\lambda \, \xi(\lambda), \\
|\gamma|
\int_{-\infty}^\infty dt \, \mu_t^{\text{ac}}(\Delta)
&=\int_{\Delta \cap \cA}d\lambda \, \xi(\lambda),
\end{align}
with
\begin{equation}
\xi(\lambda)=1+\frac{1}{\pi} \Im \bigg (\log\bigg(\frac{
\pi m_0(\lambda)+2\sqrt{|\det(X)|}}{\pi
m_0(\lambda)-2\sqrt{|\det(X)|}}\bigg)\bigg) \text{ for
a.e.~$\lambda\in\bbR$} 
\end{equation}
and $m_0$ given by \eqref{3.37}. Moreover, 
\begin{align}
|\gamma|\int_{-\infty}^\infty dt \, \mu_t^{\text{sc}}(\Delta)&=
|\Delta\cap\cS(M_0) \cap \cR(M_0)|
\intertext{and} 
|\gamma|\int_{-\infty}^\infty dt \, \mu_t^{\text{pp}}(\Delta)&=
|\Delta\cap\cP(M_0) \cap\cR(M_0)|,
\end{align}
where
\begin{equation}
\cR(M_0)=\Big \{\lambda\in\bbR \,\,\Big |\,\,  \lim_{\varepsilon
\downarrow 0}M_0(\lambda+i\varepsilon)\in \bbR \text{
and } |m_0(\lambda)|< (2/\pi)\sqrt{|\det(X)|}\Big \}. 
\end{equation}
\end{theorem}

\begin{remark}
Theorem \ref{t3.6} shows, in particular, the universality
of  the averaging on the whole parameter space in the case of cyclic
groups $g_t$ (associated with the
one-parameter subgroups $e^{tX}$ of $SL_2(\bbR)$), where $\det(X)>0$,
and in their 
limiting cases corresponding to $\det(X)=0$. However, the theorem
also shows that averaging in the case of hyperbolic one-parameter 
subgroups $g_t$, with $\det(X)<0$, depends on the initial Herglotz
function $M_0$.
\end{remark}
\begin{remark}
An analogous result concerning the decomposition of Lebesgue measure
$|\cdot |$ restricted to $\cA$ and $|\cdot |$ restricted to
$\bbR\backslash\cA$ into integrals of the measures $\mu_t^{ac}$ and
$\mu_t^{sing}=\mu_t^{sc}+\mu_t^{sing}$ on the unit circle first 
appeared in \cite{Al87}. In the case of self-adjoint rank-one
perturbations of self-adjoint operators (which is a special case of
\eqref{3.108} as observed in Remark \ref{r2.4}), \eqref{3.108} appeared
in \cite{DSS94}. 
\end{remark}

\section{Spectral averaging and Hausdorff measures} \lb{s4}

Lebesgue's decomposition of measures \eqref{leb} is a particular case of a
more general result in the theory of  decomposing measures with respect to
Hausdorff measures. This result states, in particular, that for each
$\kappa\in [0,1]$, a  Borel measure $\mu$ can be decomposed uniquely as 
\begin{equation}
\mu=\mu^{\kappa\text{-}c} +\mu^{\kappa\text{-}s}, \lb{4.1}
\end{equation}
where $\mu^{\kappa\text{-}c}$ is $\kappa$-continuous with respect to the
$\kappa$-dimensional Hausdorff measure $h^\kappa$ (i.e., 
$\mu^{\kappa\text{-}c}$
gives zero weight to sets with zero $\kappa$-dimensional Hausdorff
measure $h^\kappa$) and $\mu^{\kappa\text{-}s}$ is $\kappa$-singular with respect
to the $\kappa$-dimensional Hausdorff measure
(i.e., $\mu^{\kappa\text{-}s}$ is supported on a set with of zero 
$\kappa$-dimensional Hausdorff measure $h^\kappa$). For more details
on the decomposition \eqref{4.1} we refer to \cite{La96}, \cite{RT59},
\cite{RT63}.

\medskip

We recall that the $\kappa$-dimensional Hausdorff (outer) measure
$h^\kappa$, $\kappa\in [0,1]$ of a set $S\subset \bbR$ 
is defined as 
\begin{equation}
h^\kappa(S)=\lim_{\delta\downarrow 0}
\inf_{\delta-\text{covers}}\sum_{n\in \bbN} |I_n(\delta)|^\kappa, 
\end{equation}
where the infimum is taken over countable collections of intervals
$\{I_n(\delta)\}_{n\in \bbN}$, the $\delta$-covers, such that 
\begin{equation}
S\subset\bigcup_{n\in \bbN}I_n(\delta)\quad \text{ and }
|I_n(\delta)|<\delta \text{ for all } n\in\bbN. 
\end{equation}
We also recall that the Hausdorff dimension of a set $S$  is defined by
\begin{equation}
\dim_H(S)=\inf \{\kappa\in [0,1] \,\,\vert \,\, h^\kappa(S)=0\}. 
\end{equation}

The goal of this section is to obtain partial results concerning
spectral averaging of the $\kappa$-continuous part
$\mu_t^{\text{sc},\kappa\text{-}c}$ with respect to $h^\kappa$,
$\kappa\in (0,1)$ of the singular continuous part $\mu_t^{sc}$
(with respect to Lebesgue measure) \eqref{leb} of the measure 
$\mu_t$ associated with the family of Herglotz functions $M_t=g_t(M_0)$, 
\begin{equation}
\mu_t^{sc}=\mu_t^{\text{sc},\kappa\text{-}c}
+\mu_t^{\text{sc},\kappa\text{-}s}, \quad t\in \bbR, 
\end{equation}
where $g_t$ is a one-parameter group of automorphisms of $\Aut(\bbC_+)$.

We introduce the following hypothesis.

\begin{hypothesis}\label{fractal}
 Let $M_0$ be a Herglotz function of the type \eqref{hrhr}, $\kappa\in
(0,1)$, 
\begin{equation}
\cS_\kappa(M_0)=\Big \{\lambda\in \bbR\,\,\Big|\,\, \liminf_{\varepsilon
\downarrow 0}\varepsilon^{\kappa-1}\Im(M_0(\lambda+i\varepsilon))\in(0,
\infty) \Big\}, 
\end{equation}
and assume that the set $\cA_\kappa(M_0)$, defined by
\begin{equation}
\cA_\kappa(M_0)=\bigcup_{\kappa'\in [\kappa,1)}\cS_{\kappa'}(M_0), 
\end{equation}
is a Borel set of positive Lebesgue measure.
\end{hypothesis}

We note that by Hypothesis \ref{fractal},
\begin{equation}\label{sinsin}
\cA_\kappa(M_0)\subseteq \cS(M_0),
\end{equation}
where $\cS(M_0)$ is the invariant set \eqref{invsin} associated with the
Herglotz function $M_0$.

\begin{lemma}\label{safractal}
Assume Hypothesis \ref{fractal} and the hypotheses of Theorem
\ref{t3.6}.
Let
\begin{equation}
\mu_t^{sc}=\mu_t^{\text{sc},\kappa\text{-}c}
+\mu_t^{\text{sc},\kappa\text{-}s}, \quad t\in \bbR
\end{equation}
be the decomposition of the measure $\mu_t^{sc}$ \eqref{leb} such that  
$\mu^{\text{sc},\kappa\text{-}c}$ is $\kappa$-continuous and
$\mu^{\text{sc},\kappa\text{-}s}$ is $\kappa$-singular $($with respect to
the $\kappa$-dimensional Hausdorff measure $h^\kappa$$)$. Then, for any
bounded Borel set $\Delta\subset \bbR$  and
$0\ne|t|<\pi{\big/}\big(2\sqrt{\det(X)}\big)$ in case I, and $0\ne t\in
\bbR$ in cases II and III, 
 \begin{equation}\label{dost}
\mu_t^{sc}(\Delta\cap\cA_\kappa(M_0))
=\mu_t^{\text{sc},\kappa\text{-}c}(\Delta\cap
\cA_\kappa(M_0)).
\end{equation}
\end{lemma}
\begin{proof}
We note that 
\begin{equation}
0<\liminf_{\varepsilon \downarrow 0}
\varepsilon^{\kappa-1}\Im(M_0(\lambda+i\varepsilon)) \text{ (possibly
equal to $+\infty$) for } \lambda\in \cA_\kappa(M_0). 
\end{equation}
Using the estimate
\begin{equation}
\Im(M_t(z))=\frac{\Im(M_0(z))}{|c_tM_0(z)+d_t|^2}\le\frac{1}{c_t^2
\Im(M_0(z))}, \quad
z\in \bbC_+
\end{equation}
(we recall that $c_t\ne 0$ and $d_t\in\bbR$ by hypothesis), one infers
\begin{equation}\label{meas}
0\le \limsup_{\varepsilon \downarrow 0}
\varepsilon^{1-\kappa}\Im(M_t(\lambda+i\varepsilon))<\infty, \quad 
\lambda\in \cA_\kappa(M_0). 
\end{equation}
It is known (cf.\ \cite[Lemma\ 3.2]{DJLS96}) that \eqref{meas} implies 
\begin{equation}
0\le\limsup_{\delta \downarrow 0}
\frac{\mu_t(\lambda-\delta, \lambda+\delta)}{\delta^\kappa}
<\infty, \quad \lambda\in \cA_\kappa(M_0). 
\end{equation}
Hence, by a result of Rogers-Taylor \cite{RT59}, \cite{RT63} (also see
\cite[ Theorem 2.1]{DJLS96}) the measure 
$\mu_t \upharpoonright \cA_\kappa(M_0)$ is a $\kappa$-continuous 
measure, which proves \eqref{dost}, since 
 \begin{equation}
 \mu_t \upharpoonright \cA_\kappa(M_0)=\mu_t^{sc} \upharpoonright
\cA_\kappa(M_0)
\end{equation}
by \eqref{sinsin}.
\end{proof}

\begin{remark}
In general, we can neither state that $\cA_\kappa(M_0)$ is a Borel set
$($cf.\ Hypothesis \ref{fractal}$)$, nor that
\begin{equation}
\mu_t^{sc}(\Delta)
=\mu_t^{\text{sc},\kappa\text{-}c}(\Delta\cap \cA_\kappa(M_0)),
\quad t\ne 0. 
\end{equation}
\end{remark}

It was pointed out to us by Barry Simon that a different but not
unrelated discussion of singular continuous measures for continuous and
discrete half-line Schr\"odinger operators, based on asymptotic
behavior of solutions, was recently provided in \cite{KLS01}
(following a previous result in \cite{JL99}).  

As a corollary we get the following result. 

\begin{corollary} \lb{c4.4}
Assume the hypotheses of Lemma \ref{safractal}. Then for any bounded
Borel
set $\Delta\subset \bbR$ the following hold. \\
$(i)$ If $\det(X)>0$ then,
\begin{equation}
|\gamma|\int_{-\pi{\big/}\big(2\sqrt{\det(X)}\big)}^{\pi{\big/}
\big(2\sqrt{\det(X)}\big)} dt
\, \mu_t^{\text{sc},\kappa\text{-}c}(\Delta\cap \cA_\kappa(M_0)) 
=|\Delta\cap \cA_\kappa(M_0)|. \lb{4.17}
\end{equation}
$(ii)$ If $\det(X)=0$ then,
\begin{equation}
|\gamma|\int_{-\infty}^\infty dt
\, \mu_t^{\text{sc},\kappa\text{-}c}(\Delta\cap \cA_\kappa(M_0)) =
|\Delta \cap \cA_\kappa(M_0)|. \lb{4.18}
\end{equation}
$(iii)$ If $\det(X)<0$ then, 
\begin{equation}
|\gamma|\int_{-\infty}^\infty dt \, 
\mu_t^{\text{sc},\kappa\text{-}c}(\Delta\cap
\cA_\kappa(M_0))=
|\Delta\cap \cA_\kappa(M_0)\cap \cR(M_0)|, \lb{4.19}
\end{equation}
where
\begin{align}
\cR(M_0)&=\Big \{\lambda\in\bbR \,\Big|\,  \lim_{\varepsilon \downarrow
0}M_0(\lambda+i\varepsilon)\in \bbR \text{ and } \no \\ 
& \hspace*{.7cm} 
\big|\lim_{\varepsilon \downarrow 0}(\gamma M_0(\lambda+i\varepsilon)
-\beta)\big| < (2/\pi)\sqrt{|\det(X)|}\Big \}. \lb{4.20}
\end{align}
\end{corollary}

Even though Corollary \ref{c4.4}
appears to be a new result, it cannot be considered a complete analog
of
\eqref{3.61}, since first of all we have no results for the singular
part
$\mu_t^{\text{sc},\kappa-\text{s}}$, and secondly, we were not able to
remove the set $\cA_\kappa(M_0)$ on the left-hand sides of
\eqref{4.17}--\eqref{4.19}. We hope our present attempt will encourage
future work in this direction. \\

\noindent {\bf Acknowledgments.} We are indebted to Vadim
Kostrykin, Yuri Latushkin, and Barry Simon for stimulating discussions
and hints  to pertinent literature and grateful to the referee for
a careful reading of the manuscript. \\


\end{document}